\documentclass[11pt]{amsart}

\usepackage{latexsym}
\usepackage{psfrag}
\usepackage{amsmath}
\usepackage{amssymb}
\usepackage{epsfig}
\usepackage{amsfonts}
\usepackage{amscd}

\def\RR{Rogers-Ramanujan }
\def\al{\alpha}
\def\be{\beta}

\def\la{\lambda}
\def\si{\sigma}
\def\pr{\prime}

\def\zz{\mathbb Z}

\newcommand{\qbinom}[2]{\genfrac{[}{]}{0pt}{}{#1}{#2}_q}

\newtheorem{theorem}{Theorem}[section]
\newtheorem{lemma}[theorem]{Lemma}
\newtheorem{obs}[theorem]{Observation}

\newtheorem{cor}[theorem]{Corollary}
\newtheorem{prop}[theorem]{Proposition}

\newtheorem{defn}[theorem]{Definition}

\newtheorem{alg}[theorem]{Procedure}

\begin{document}%%%%%%%%%%%%%%%%%%%%%%%%%%%%%%%%%%%%%%%%%%%%%%%%%%%%%%%%
%%%%%%%%%%%%%%%%%%%%%%%%%%%%%%%%%%%%%%%%%%%%%%%%%%%%%%%%%%%%%%%%%%%%%%%%

\title[Dyson's new symmetry and \\ generalized Rogers-Ramanujan identities]%
    {Dyson's new symmetry and \\ generalized Rogers-Ramanujan identities}
\author[Cilanne Boulet ${\hspace{-.95ex}}^\ast$]{Cilanne Boulet ${\hspace{-.95ex}}^\ast$}
% \date{November 29, 2005}


\thanks{${\hspace{-.95ex}}^\ast$Department of Mathematics, Cornell University, Ithaca, NY
14853. 
\ Email: \texttt{cilanne@math.cornell.edu}}

\begin{abstract}  
\noindent  
We present a generalization, which we call $(k,m)$-rank, of Dyson's notion
of rank to integer partitions with $k$ successive Durfee rectangles and give two
combinatorial symmetries associated with this new definition.  We prove
these symmetries bijectively.  Using the two symmetries we give a new
combinatorial proof of generalized Roger-Ramanujan identities.  We also
describe the relationship between $(k,m)$-rank and Garvan's $k$-rank.  
\end{abstract}

\maketitle

%%%%%%%%%%%%%%%%%%%%%%%%%%%%%%%%%%%%%%%%%%%%%%%%%%%%%%%%%%%
%%%%                                                   %%%%
%%%%                   Introduction                    %%%%
%%%%                                                   %%%%
%%%%%%%%%%%%%%%%%%%%%%%%%%%%%%%%%%%%%%%%%%%%%%%%%%%%%%%%%%%

\section{Introduction} \label{sec:intro}

First discovered by Rogers~\cite{Rogers}  in 1894, the \RR identities,
$$
\sum_{n=0}^\infty \frac{q^{n^2}}{(1-q)(1-q^2)\cdots (1-q^n)} \ = \
\prod_{n=0}^\infty \frac{1}{(1-q^{5n+1})(1-q^{5n+4})} \, 
$$
and
$$
\sum_{n=0}^\infty \frac{q^{n^2+n}}{(1-q)(1-q^2)\cdots (1-q^n)} \ = \
\prod_{n=0}^\infty \frac{1}{(1-q^{5n+2})(1-q^{5n+3})} \, ,
$$
are among the most intriguing partition identities.  

The goal of this paper is to present a new combinatorial proof of the
following generalization (which is due
to Andrews~\cite{Andr:odd}) of the first \RR identity, for $k \geq 1$:
\begin{equation}
\label{eqn:genRR}
\sum_{n_1 = 0}^\infty \cdots \sum_{n_{k-1} = 0}^\infty
\frac{q^{N_1^2+N_2^2+\dots+N_{k-1}^2}}{(q)_{n_1}(q)_{n_2}\dots(q)_{n_{k-1}}} \  =
\hspace{-3mm} \prod_{
\scriptsize{\begin{array}{c}
n=1 \\
n \not\equiv 0, \pm k~(\mathrm{mod~} 2k+1)
\end{array}}}^{\infty} \hspace{-4mm}\frac{1}{1-q^n}\hspace{-5mm}
\end{equation}
where~$N_j = n_j + n_{j+1} + \cdots + n_{k-1}$.  We use the standard $q$-series notation and let 
$(q)_\infty = \prod_{i=1}^\infty (1-q^i)$ and $(q)_n = \prod_{i=1}^n (1-q^i)$.

Instead of attacking this identity directly, we will use two bijections to prove the following family of identities, which we
call the generalized Schur identities:
\begin{equation}
\label{eqn:genSchur}
\sum_{n_1 = 0}^\infty \cdots \sum_{n_{k-1} = 0}^\infty \frac{q^{N_1^2+N_2^2+\dots+N_{k-1}^2}}{(q)_{n_1}(q)_{n_2}\dots(q)_{n_{k-1}}} =
\frac{1}{(q)_\infty} \sum_{j = -\infty}^{\infty} (-1)^j q^{\frac{j(j+1)(2k+1)}{2}-kj} \, ,
\end{equation}
with~$N_j = n_j + n_{j+1} + \cdots + n_{k-1}$.  
By
using Jacobi's triple product identity,
$$
\sum_{j=-\infty}^\infty z^j t^{\frac{j(j+1)}{2}}  \,
= \, \prod_{i=1}^\infty (1+zt^i) \,
\prod_{j=0}^\infty (1+z^{-1}t^{j}) \,
\prod_{i=1}^\infty (1-t^i) \, ,
$$
which specializes to
$$
\sum_{j=-\infty}^\infty (-1)^j q^{\frac{j(j+1)(2k+1)}{2}-kj}  \,  =
\hspace{-3mm} \prod_{
\scriptsize{\begin{array}{c}
n=1 \\
n \equiv 0, \pm k~(\mathrm{mod~} 2k+1)
\end{array}}}^{\infty} \hspace{-7mm}{1-q^n} \, 
$$
when we let $t = q^{2k+1}$ and $z = -q^{-k}$,
we see that~(\ref{eqn:genRR}) and~(\ref{eqn:genSchur}) are equivalent.  This application
of Jacobi's triple product identity is a standard first step in \RR proofs and in particular
was used by Schur~\cite{Schur} in his combinatorial proof of the original \RR identities. 
We note that the
Jacobi triple product identity has a combinatorial proof
due to Sylvester (see~\cite{Pak:survey,Wright:Jacobi}) and 
so its application does not change the combinatorial nature of our proof.   

Before presenting our proof of the generalized \RR identities~(\ref{eqn:genRR}), we must outline our notation and review two important
ideas.  The first is  Andrews' notion of successive Durfee squares which gives a combinatorial interpretation
to the left hand side of~(\ref{eqn:genRR}) and~(\ref{eqn:genSchur}).  The second is Dyson's proof of Euler's pentagonal
number theorem based on his definition of rank.

\subsection{Notation}
We
begin by giving the basic definitions that we will need. This section is meant simply to
familiarize the reader with the notation that will be used, rather than provide an
introduction to the subject.  For such an introduction we recommend~\cite{Andr:book,
Pak:survey}.  

A {\em partition}~$\la$ is a sequence of integers~$(\la_1, \la_2, ..., \la_\ell)$ such
that~$\la_1 \geq  \la_2 \geq ... \geq  \la_\ell > 0$.   As a convention, we will say
that $\la_j = 0$ for $j > \ell$.   We call each~$\la_i$ a {\em
part} of~$\la$.   We say that~$\la$ is a partition of~$n$, denoted~$\la \vdash n$
or~$|\la| = n$, if~$\sum \la_i = n$.  Let $\mathcal{P}_n$ denote the set of partitions
of $n$ and let $p(n) = |\mathcal{P}_n|$.  Also, let $\mathcal{P} = \cup_n \mathcal{P}_n$
denote the set of all partitions.  

We let~$\ell(\la) = \ell$ denote the {\em number of non-zero parts} of~$\la$.  In general, we will let
number of part mean number of non-zero parts.  Let~$f(\la) =
\la_1$ denote the {\em largest part} of~$\la$, and let~$e(\la) = \la_{\ell(\la)}$ denote
the {\em smallest (non-zero) part} of $\la$.   

\begin{figure}[hbt]
\begin{center}
\psfrag{1}{$\la$}
\psfrag{2}{$\la^\pr$}
\epsfig{file=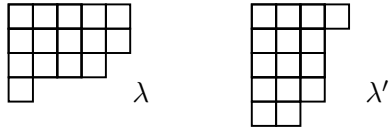,height=1.666cm}
\end{center}
\caption{Partition $\la = (5,5,4,1)$ and conjugate partition
$\la^\pr = (4,3,3,3,2)$.}
\label{fig:intro:part}
\end{figure}

To every partition we associate a {\em Young diagram} as in Figure~\ref{fig:intro:part}.  

The {\em conjugate}~$\la^\pr$ of a partition~$\la$ is obtained by reflection across the main
diagonal (again see Figure~\ref{fig:intro:part}).  Alternatively,~$\la^\pr$ may be defined as
follows: $\la^\pr = (\la^\pr_1, \la^\pr_2, ..., \la^\pr_{f(\la)})$ where 
$\la^\pr_i = |\{j:\la_j \geq i\}|$ 
is the number of parts of~$\la$ which are greater than or equal to~$i$.

\subsection{Andrews' successive Durfee squares}

Andrews introduced the idea of successive Durfee squares to study his generalized Rogers-Ramanujan identities~\cite{Andr:Durfee}.  
He interpreted the left hand sides of equations~(\ref{eqn:genRR}) and~(\ref{eqn:genSchur})
as follows.

\begin{figure}[hbt]
\begin{center}
\psfrag{l}{$\la$}
\epsfig{file=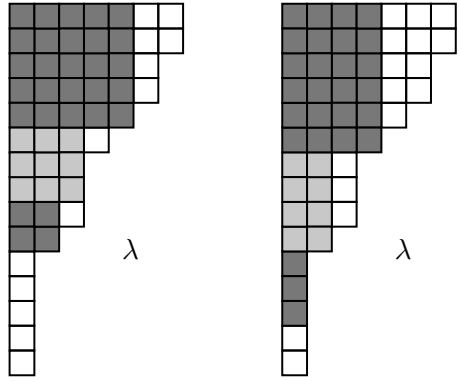,height=5cm}
\end{center}
\caption{The first three successive Durfee squares and $1$-Durfee rectangles of $\la = (7,7,6,6,5,4,3,3,3,2,1,1,1,1,1)$.  On the left we see that $\la$ has successive Durfee squares of size 5, 3, and 2.  On the right we see that $\la$ has successive $1$-Durfee rectangles of width 4, 2, and 1.  }
\label{figSuccDur}
\end{figure}

\begin{defn}
The {\em first Durfee square} of a partition~$\la$ is the largest square that fits in the upper left hand corner of the diagram of~$\la$.  The {\em second Durfee square} is the largest square that fit in the diagram of~$\la$ below the first Durfee square of~$\la$.  In general, the {\em $k$th Durfee square} is the largest square that fits below the $(k-1)$st Durfee square of~$\la$.  
\end{defn}
See Figure~\ref{figSuccDur} for an example.

Let $q_k(n)$ denote the number of partitions of $n$ with at most $k$ Durfee squares
and let $\mathcal{Q}_{k}$ denote the set of all partitions with at most $k$ Durfee squares.
Now the generating function for partitions with Durfee squares of size $N_1$, $N_2$, ..., $N_{k-1}$ and no part below the $k-1$st Durfee square is
$$
\frac{q^{N_1^2+N_2^2+\dots+N_{k-1}^2}}{(q)_{n_1}(q)_{n_2}\dots(q)_{n_{k-1}}}
$$
where $n_j = N_j - N_{j+1}$ so that~$N_j = n_j + n_{j+1} + \cdots + n_{k-1}$.  
This can be seem by a simple counting argument as is done by Andrews~\cite{Andr:Durfee}. 
Alternatively, in Appendix A of~\cite{Bou:Thesis}, we show this bijectively using the insertion
procedure which is defined in this paper.  

Therefore the generating function for partitions with at most $k-1$ Durfee squares is 
$$
1 + \sum_{n=1}^\infty q_{k-1}(n) \, q^n \, = \, \sum_{n_1 = 0}^\infty \cdots \sum_{n_{k-1} = 0}^\infty \frac{q^{N_1^2+N_2^2+\dots+N_{k-1}^2}}{(q)_{n_1}(q)_{n_2}\dots(q)_{n_{k-1}}}
$$
with ~$N_j = n_j + n_{j+1} + \cdots + n_{k-1}$
which is indeed the left hand side of~(\ref{eqn:genRR}).  

For our proof, we extend the notion of successive Durfee squares.  
\begin{defn}
For any integer $m$, define an {\em $m$-rectangle} to be a rectangle whose height exceeds its
width by exactly $m$.   We require $m$-rectangles to have non-zero height though they may have
width zero.  
\end{defn}
In particular, notice that $0$-rectangles are simply squares.  The technical detail about
zero width being allowed will be used to obtain Observation~\ref{genObs2}.

We define successive $m$-Durfee rectangles in the same manner as Andrews' successive Durfee squares.  
\begin{defn}
The {\em first $m$-Durfee rectangle} of a partition~$\la$ is the largest $m$-rectangle that
fits in the upper left hand corner of the diagram of~$\la$.  The {\em second $m$-Durfee
rectangle} is the largest $m$-rectangle that fits in the diagram of~$\la$ below the first
Durfee square of~$\la$.  In general, the {\em $k$th successive $m$-Durfee rectangle} is the
largest $m$-rectangle that fits below the $(k-1)$st Durfee square of~$\la$.
\end{defn}
Again, see Figure~\ref{figSuccDur} for an example of successive $1$-Durfee rectangles.

Note that the possibility of width zero $m$-rectangles means that, for $m > 0$, all partitions
(including the empty partition) have arbitrarily many successive $m$-Durfee rectangles.  
In this case, the
Durfee rectangles extend below the non-zero parts of the partition.

\subsection{Dyson's proof of Euler's pentagonal number theorem}
The primary inspiration for the \RR proof in this paper is Dyson's proof of Euler's pentagonal number 
theorem,
$$
1 =
\frac{1}{(q)_\infty} \sum_{j = -\infty}^{\infty} (-1)^j q^{\frac{j(3j-1)}{2}} \, ,
$$
based on his definition of rank~\cite{Dyson:new} (see also~\cite{Dyson:walk}).
Note that this identity is the case $k=1$ of the generalised Schur identities~(\ref{eqn:genSchur}).

\begin{defn}[Dyson~\cite{Dyson:guesses}]
The rank of a partition $\la$ is 
$$
r(\la) = f(\la) - \ell(\la) .
$$
\end{defn}
Recall that $f(\la)$ is the length of the first part of $\la$ and $\ell(\la)$ is the number of parts of~$\la$.  

\begin{figure}[hbt]
\begin{center}
\psfrag{1}{\, $\la$}
\psfrag{2}{\hspace{-2.4mm}$d_{-2}(\la)$}
\psfrag{3}{$d_1(\la)$}
\epsfig{file=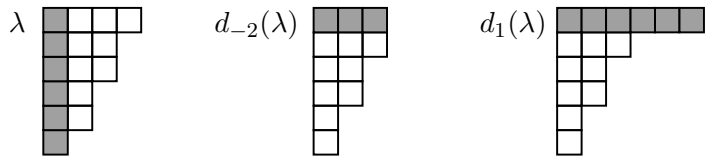,height=2cm}
\end{center}
\caption{Applying Dyson's map to $\la = (4,3,3,2,2,1)$ with $r(\la) = 4-6 = -2$ gives $d_{-2}(\la) = (3,3,2,2,1,1)$ and $d_{1}(\la) = (6,3,2,2,1,1)$.}
\label{figDysEg}
\end{figure}

Dyson's proof of Euler's pentagonal number theorem can be presented as follows.  (In addition to his 
papers, see~\cite{BG:obs} and~\cite{Pak:missed} for additional descriptions.)   
Let $h(n,r)$, $h(n,\leq r)$, and $h(n,\geq r)$ denote the number of partitions of $n$ with rank equal to $r$, less than or equal to $r$, and greater than or equal to $r$ respectively.  Clearly, for $n > 0$, we observe that 
$$
p(n) = h(n, \leq r) + h(n, \geq r+1)
$$
and Dyson noticed two symmetries,
$$
\aligned
h(n,r) & = h(n,-r) \, \, \, \text{and}\\
h(n, \leq r) & = h(n + r - 1, \geq r -2) \, .
\endaligned
$$
The first of these symmetries is a simple consequence of conjugation.  The second symmetry, the ``new symmetry'' from the title of~\cite{Dyson:new}, follows from a bijection, $d_r$, which we call Dyson's map.  Dyson's map $d_r$ takes a partition~$\la$ of~$n$ with $r(\la) \leq r$ and returns a partition~$\mu$ of~$n+r-1$ with $r(\mu) \geq r-2$ by removing the first column of~$\la$, which has~$\ell(\la)$ squares, and adding a part of size $\ell(\la) +r -1$.  This new part will be the first row of~$\mu$.  See figure~\ref{figDysEg} for an example.

Let 
$$
\aligned
H_{\leq r}(q) &  = \sum_{n=1}^\infty h(n,\leq r) q^n \, \text{ and}\\
H_{\geq r}(q) &  = \sum_{n=1}^\infty h(n,\geq r) q^n
\endaligned
$$
be the generating functions for partitions with rank at most $r$ and at least $r$.  Then
$$
H_{\leq r}(q) = q^{1-r} H_{\geq r-2}(q) = q^{1-r}\left(\frac{1}{(q)_\infty} - H_{\leq r-3}(q)\right)
$$
where the first equality follows from Dyson's new symmetry and the second equality follows from the observation.  
Applying this equation repeatedly gives
$$
\aligned
H_{\leq r}(q) & = q^{1-r}\left(\frac{1}{(q)_\infty} - H_{\leq r-3}(q)\right) \\[2mm]
              & = q^{1-r}\left(\frac{1}{(q)_\infty}\right) - q^{5-2r}\left(\frac{1}{(q)_\infty} - H_{\leq r-6}(q)\right) \\[2mm]
              & = q^{1-r}\left(\frac{1}{(q)_\infty}\right) - q^{5-2r}\left(\frac{1}{(q)_\infty}\right) + q^{12-3r}\left(\frac{1}{(q)_\infty} - H_{\leq r-9}(q)\right) \\[2mm]
              & \qquad \vdots \\[2mm]
              & = \frac{1}{(q)_\infty} \sum_{j=1}^\infty (-1)^{j-1} q^{\frac{j(3j-1)}{2}-jr} \, .
\endaligned
$$
Finally, the first symmetry (conjugation) gives us
$$
\frac{1}{(q)_\infty} = 1 + H_{\leq 0}(q) + H_{\geq 1}(q) = 1 + H_{\leq 0}(q) + H_{\leq -1}(q)
$$
and substituting gives Euler's pentagonal number theorem.  

\subsection{Outline of our proof}

Roughly speaking, our proof of the generalized Schur identities is a Dyson-style proof with a modified Dyson's rank.    
In section~\ref{secTools}, we develop the basic tools needed for our proof, selection and insertion.  
In section~\ref{secDef}, we generalized Dyson's rank to partitions with $k$ successive $m$-Durfee rectangles.  This new rank will be called
$(k,m)$-rank.  The definition will use the selection procedure from section~\ref{secTools}.
Similarly to the case of Dyson's rank, $(k,m)$-rank will satisfy two symmetries.  We prove these in section~\ref{secSym} by two bijections
that are build using using selection and insertion.  The first bijection generalizes conjugation and the second bijection generalizes the
map $d_r$ which corresponds to Dyson's new symmetry.  In section~\ref{secAlg}, we use the same algebraic
manipulations used to deduce
Euler's pentagonal number theorem to deduce the generalized Schur identities.  
We conclude by explaining how our work relates to the work of others and by mentioning a problem that is
still open.

%%%%%%%%%%%%%%%%%%%%%%%%%%%%%%%%%%%%%%%%%%%%%%%%%%%%%%%%%%%
%%%%                                                   %%%%
%%%%       Basic Tools: Selection and Insertions       %%%%
%%%%                                                   %%%%
%%%%%%%%%%%%%%%%%%%%%%%%%%%%%%%%%%%%%%%%%%%%%%%%%%%%%%%%%%%

\section{Basic Tools: Selection and Insertions}
\label{secTools}

In this section, we develop the basic tools that will be needed for our proof
of the generalized Schur identities~(\ref{eqn:genSchur}).  We define two procedures which we call {\em selection} and {\em
insertion} and we will say precisely in what sense they are inverses of each other.  

\subsection{Selection of parts from a sequence of partitions $\la^1, \la^2, ..., \la^k$}

The first procedure, selection, has as input a sequence of partitions
and as output one part from each of these partitions.  

\begin{alg}
Given a sequence of $k-1$ nonnegative integers,
$$
p_2, p_3, ..., p_k,
$$
and $k$ partitions,
$$
\la^{1}, \la^{2}, ..., \la^{k},
$$
such that
$$
f(\la^{2}) \leq p_2, f(\la^{3}) \leq p_3, ..., f(\la^{k}) \leq p_k \, ,
$$
we select one row from each partition as follows:
\begin{itemize}
\item
select the first (that is, the largest) part of $\la^k$,
\item
suppose we have selected the~$j$th part,~${\la^i}_j$, of~${\la^i}$, then select the 
$(j+p_i-{\la^i}_j)$th part 
of~$\la^{i-1}$. 
\end{itemize}
\end{alg}

One way to think of the selection of the $(j +p_i - \la^i_j)$th part of $\la^{i-1}$ is to this
that we are selecting the row of $\la^{i-1}$ that is $p_i - \la^i_j$ lower in the Young diagram 
than the row selected in $\la^i$. 
The number $p_i - \la^i_j$ can be thought of as the number of ``missing'' boxes in the $j$th row
of $\la^i$ since $\la^i$ is restricted to having parts of size at most $p_i$.

\begin{defn}
Let $A(\la^{1}, \la^{2}, ..., \la^{k}; p_2, p_3, ..., p_k)$ be the sum of the selected parts.  
\end{defn}
When $p_2, p_3, ..., p_k$ are clear from the context, we will write $A(\la^1, \la^2, ...,
\la^k)$.

\begin{figure}[hbt]
\begin{center}
\psfrag{1}{$\la^1$}
\psfrag{2}{$\la^2$}
\psfrag{3}{$\la^3$}
\psfrag{4}{$\la^4$}
\psfrag{5}{$\mu^1$}
\psfrag{6}{$\mu^2$}
\psfrag{7}{$\mu^3$}
\psfrag{8}{$\mu^4$}
\psfrag{9}{$\mu^5$}
\epsfig{file=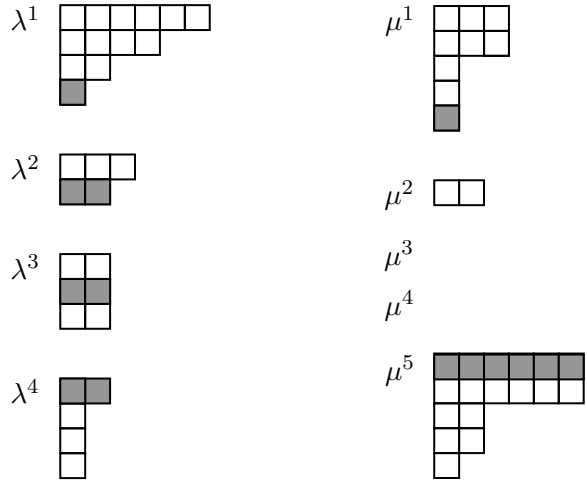,height=6.333cm}
\end{center}
\caption{Selection of rows from $\la^1$, $\la^2$, $\la^3$, and $\la^4$ with $p_2 = 4$, $p_3 = 2$,
and $p_4 = 3$ to get $A(\la^1, \la^2, \la^3, \la^4) = 7$. Selection of rows from $\mu^1$,
$\mu^2$, $\mu^3$, $\mu^4$, and $\mu^5$ with $p_2 = 2$, $p_3 = 0$, $p_4 = 2$, and $p_5 = 6$ to get
$A(\mu^1, \mu^2, \mu^3, \mu^4, \mu^5) = 7$.   Selected parts are shown in grey.}

\label{figSeqSelect}
\end{figure}

See Figure~\ref{figSeqSelect} for examples of this selection procedure.  On the left hand side,
we have $p_2 = 4$, $p_3 = 2$, and $p_4 = 3$.  We select the first part of $\la^4$.  Then we
select the $1 +(p_4-\la^4_1) = 1 +(3-2) = 2$nd part from $\la^3$, the~$2 + (p_3+\la^3_2)=2 +(2-2) = 2$nd part from $\la^2$, and the
$2+(p_2+\la^2_2) = 2+(4-2) = 4$th part from $\la^1$.  This gives $A(\la^1, \la^2, \la^3, \la^4) = 1+2+2+2 = 7$.  

On the right hand side, we have $p_2 = 2$, $p_3 = 0$, $p_4 = 2$, and $p_5 = 6$.  We select the
first part of $\mu^5$.  Then we select the $1 +(6-6)= 1$st part from $\mu^4$, the $1+(2-0) =3$th
part from~$\mu^3$, the $3+(0-0)= 3$th part from $\mu^2$, and the $3+(2-0) =  5$th part from
$\mu^1$. This gives $A(\mu^1, \mu^2, \mu^3, \mu^4, \mu^5) = 1 +0+0+0+6 = 7$.  

%\subsection{Selection and removal of parts from a sequence of partitions $\la^1, \la^2, ..., \la^k$}
It will be useful to establish notation for the selection of parts from a sequence of partitions and the
removal of those parts.  

\begin{defn}
Given a sequence of $k-1$ nonnegative integers,
$$
p_2, p_3, ..., p_k,
$$
and $k$ partitions,
$$
\la^{1}, \la^{2}, ..., \la^{k},
$$
such that
$$
f(\la^{2}) \leq p_2, f(\la^{3}) \leq p_3, ..., f(\la^{k}) \leq p_k \, ,
$$
\begin{itemize}
\item
let $\psi_1 = {\psi_1}_{\{p_2, ..., p_k\}}(\la^1, \la^2, ..., \la^k) = 
A(\la^1, \la^2, ..., \la^k)$ and 
\item
let $\psi_2 = {\psi_2}_{\{p_2, ..., p_k\}}(\la^1, \la^2, ..., \la^k) 
= (\mu^1, \mu^2, ..., \mu^k)$ 
where $\mu^1, \mu^2, ..., \mu^k$ are found by removing the parts of $\la^1, \la^2, ..., \la^k$ 
selected while calculating $\psi_1$.  
\end{itemize}
Let $\psi_{\{p_2, ..., p_k\}}(\la^1, \la^2, ..., \la^k) =(\psi_1; \psi_2)$.  
\end{defn}

When $\{p_2, ..., p_k\}$ are clear from context, we will write $\psi_1(\la^1, \la^2, ..., \la^k)$,
$\psi_2(\la^1, \la^2, ..., \la^k)$, and $\psi(\la^1, \la^2, ..., \la^k)$.

\subsection{Insertion into a sequence of partitions $\la^1, \la^2, ..., \la^k$}
Based on this definition of selection from a sequence of partitions, we can define an insertion
algorithm on which our two symmetries are based.  The following proposition describes the result
of insertion.  We will give a procedure for insertion after the proof of this proposition.   

\begin{prop}
\label{prop:insertion}
Given a sequence of $k-1$ nonnegative integers
$$
p_2, p_3, ..., p_k,  
$$
$k$ partitions
$$
\la^{1}, \la^{2}, ..., \la^{k}
$$
with $|\la^{1}| + | \la^{2} | + ... +| \la^{k} | = n$, \\
such that
$$
f(\la^{2}) \leq p_2, f(\la^{3}) \leq p_3, ..., f(\la^k) \leq p_k \, ,
$$
and an integer $a \geq A(\la^1, \la^2, ..., \la^k; p_2, p_3, ..., p_k)$, \\
there exists a {\em unique} sequence of $k$ partitions,
$$
\mu^1, \mu^2, ..., \mu^k,
$$
obtained by inserting one (possibly empty) part into each of the original partitions,
$$
\la^1, \la^2, ..., \la^k,
$$
such that
\begin{enumerate}
\item
$|\mu^1|+ |\mu^2| + ... + |\mu^k|= n+a$,
\item
$f(\mu^2) \leq p_2, f(\mu^3) \leq p_3, ..., f(\mu^k) \leq p_k$,
\item
$A(\mu^1, \mu^2, ..., \mu^k; p_2, p_3, ..., p_k) = a$.
\end{enumerate}
Moreover, the inserted parts have the same length as those which are selected when calculating $A(\mu^1, \mu^2, ..., \mu^k; p_2, p_3, ..., p_k)$.
\end{prop}

We will prove this proposition by induction on $a$.  The two following lemmas are the required base case and inductive step.  

\begin{lemma}
\label{lemma:lem1}
Proposition~\ref{prop:insertion} (without uniqueness) is true for $a = A(\la^1, \la^2, ..., \la^k; p_2, p_3, ..., p_k).$
\end{lemma}

\begin{proof}
For each $\la^i$, consider the size of the part selected from that partition.  Insert an additional part in $\la^i$ of the same size as the selected part to obtain $\mu^i$.  See Figure~\ref{figInsertMin}.

\medskip
\begin{figure}[hbt]
\begin{center}
\psfrag{1}{$\mu^1$}
\psfrag{2}{$\mu^2$}
\psfrag{3}{$\mu^3$}
\psfrag{4}{$\mu^4$}
\psfrag{5}{$\la^1$}
\psfrag{6}{$\la^2$}
\psfrag{7}{$\la^3$}
\psfrag{8}{$\la^4$}
\epsfig{file=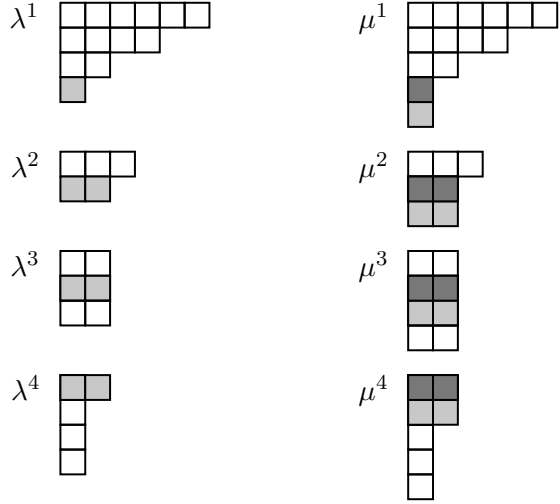,height=6.666cm}
\end{center}
\caption{With $p_2 = 4$, $p_3 = 2$, and $p_4 = 3$, inserting $7 = A(\la^1, \la^2, \la^3, \la^4)$ into $\la^1$, $\la^2$, $\la^3$, and $\la^4$ gives $\mu^1$, $\mu^2$, $\mu^3$, and $\mu^4$.}
\label{figInsertMin}
\end{figure}

We have inserted parts totaling $a = A(\la^1, \la^2, ..., \la^k; p_2, p_3, ..., p_k)$ since the sum of the selected parts of $\la^1, \la^2, ..., \la^k$ is $a$.  This implies condition (1).

Note that, for each $\la^i$, since we are inserting a part of the same size as the selected part,
it can be inserted directly above the selected row in $\la^i$.  Again since we are inserting
parts of the same size and since $p_2, p_3, ..., p_k$ remain constant,  when we select rows from
$\mu^1, \mu^2, ..., \mu^k$, we will select the rows we have just added.  Moreover, this gives
$A(\mu^1, \mu^2, ..., \mu^k; p_2, p_3, ..., p_k) = a$, condition (3).

Finally, condition (2) is satisfied since $f(\la^{i}) \leq p_i$ and the part 
selected from $\la^i$, and added to give $\mu^i$, is at most $f(\la^i)$.
\end{proof}

\begin{lemma}
\label{lemma:lem2}
If Proposition~\ref{prop:insertion} (without uniqueness) is true for $a = b$, then it is true for $a = b+1$.
\end{lemma}

\begin{proof}
Suppose $\nu^1, \nu^2, ..., \nu^k$ are the partitions obtained by inserting $b$ into $\la^1, \la^2, ..., \la^k$ as in Proposition~\ref{prop:insertion}.  To insert $b+1$ into $\la^1, \la^2, ..., \la^k$ we need to determine which partition $\la^i$ gets a part that is larger than it did when we inserted $b$ into $\la^1, \la^2, ..., \la^k$.  

If the selected part of each of $\nu^1, \nu^2, ..., \nu^k$ is the first part of that partition, 
then we let $\mu^2 = \nu^2, \mu^3 = \nu^3, ..., \mu^k = \nu^k$ and we let $\mu^1$ be $\nu^1$ 
except with first part one larger, i.e. $\mu^1_1 = \nu^1_1 + 1$ and $\mu^1_i = \nu^1_i$ for $i \geq 2$.  See Figure~\ref{figInsertTop}.

\bigskip
\begin{figure}[hbt]
\begin{center}
\psfrag{a}{$\la^1$}
\psfrag{b}{$\la^2$}
\psfrag{c}{$\la^3$}
\psfrag{1}{$\nu^1$}
\psfrag{2}{$\nu^2$}
\psfrag{3}{$\nu^3$}
\psfrag{4}{$\mu^1$}
\psfrag{5}{$\mu^2$}
\psfrag{6}{$\mu^3$}
\epsfig{file=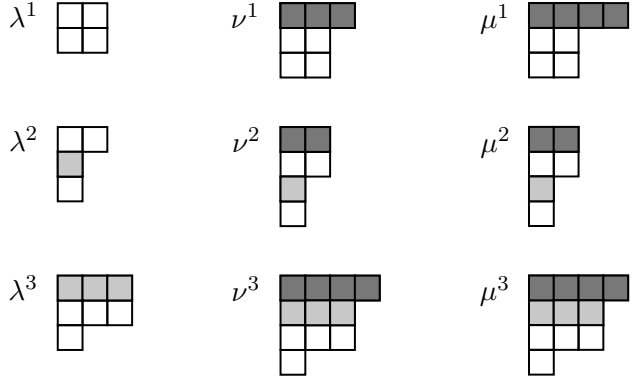,height=5cm}
\end{center}
\caption{Consider $\la^1$, $\la^2$, and $\la^3$ with $p_2 = 4$ and $p_3 = 2$ so that $A(\la^1, \la^2, \la^3) = 4$.  If $\nu^1$, $\nu^2$, and $\nu^3$ are obtained by inserting 9 into $\la^1$, $\la^2$, and $\la^3$, then $\mu^1$, $\mu^2$, and $\mu^3$ are obtained by inserting 10 into $\la^1$, $\la^2$, and $\la^3$.}
\label{figInsertTop}
\end{figure}

Otherwise consider the smallest $i$ such that the selected part of $\nu^i$ is {\em not} 
equal to the part above it or $p_i$ if it is the first part of $\nu^i$.  (Since we start by selecting
the first row of $\nu^k$, if we have not selected the first row of every $\nu^1, \nu^2, ..., \nu^k$,
there must be such an $i$.)  Add $1$ to this selected
part in $\nu^i$ to obtain $\mu^i$.  The rest of the sequence of partitions is 
defined by $\mu^j = \nu^j$.  
See Figures~\ref{figInsertMore1} and~\ref{figInsertMore2}.  

\bigskip
\begin{figure}[hbt]
\begin{center}
\psfrag{a}{$\la^1$}
\psfrag{b}{$\la^2$}
\psfrag{c}{$\la^3$}
\psfrag{d}{$\la^4$}
\psfrag{e}{$\nu^1$}
\psfrag{f}{$\nu^2$}
\psfrag{g}{$\nu^3$}
\psfrag{h}{$\nu^4$}
\psfrag{1}{$\mu^1$}
\psfrag{2}{$\mu^2$}
\psfrag{3}{$\mu^3$}
\psfrag{4}{$\mu^4$}
\psfrag{5}{$\rho^1$}
\psfrag{6}{$\rho^2$}
\psfrag{7}{$\rho^3$}
\psfrag{8}{$\rho^4$}
\epsfig{file=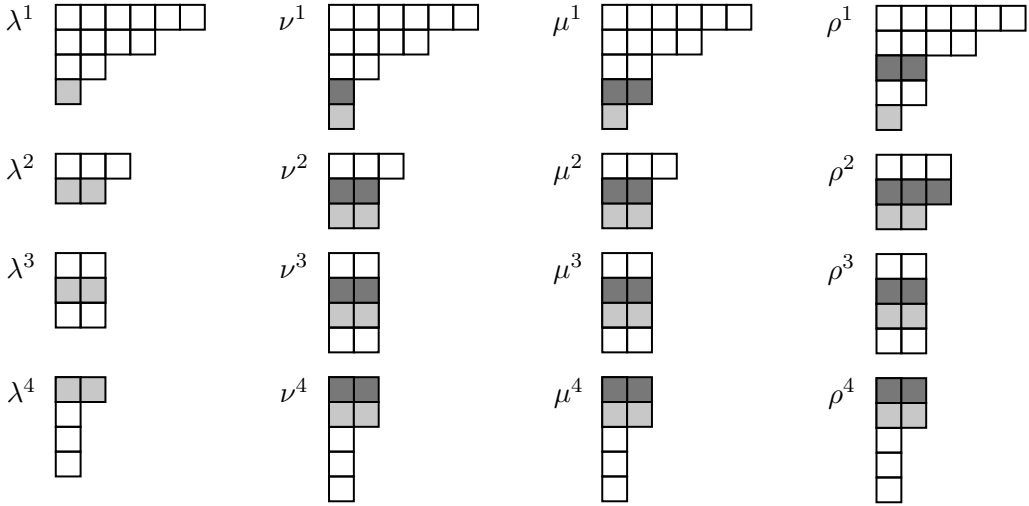,height=6.666cm}
\end{center}
\caption{Consider $\la^1$, $\la^2$, $\la^3$, and $\la^4$ with $p_2 = 4$, $p_3 = 2$, and $p_4 = 3$.  We insert 7 into $\la^1$, $\la^2$, $\la^3$, and $\la^4$ to get $\nu^1$, $\nu^2$, $\nu^3$, and $\nu^4$, 8 to get $\mu^1$, $\mu^2$, $\mu^3$, and $\mu^4$, and 9 to get $\rho^1$, $\rho^2$, $\rho^3$, and $\rho^4$. }
\label{figInsertMore1}
\end{figure}

\bigskip
\begin{figure}[hbt]
\begin{center}
\psfrag{1}{$\la^1$}
\psfrag{2}{$\la^2$}
\psfrag{3}{$\la^3$}
\psfrag{4}{$\nu^1$}
\psfrag{5}{$\nu^2$}
\psfrag{6}{$\nu^3$}
\psfrag{7}{$\mu^1$}
\psfrag{8}{$\mu^2$}
\psfrag{9}{$\mu^3$}
\psfrag{a}{$\rho^1$}
\psfrag{b}{$\rho^2$}
\psfrag{c}{$\rho^3$}
\psfrag{d}{$\si^1$}
\psfrag{e}{$\si^2$}
\psfrag{f}{$\si^3$}
\epsfig{file=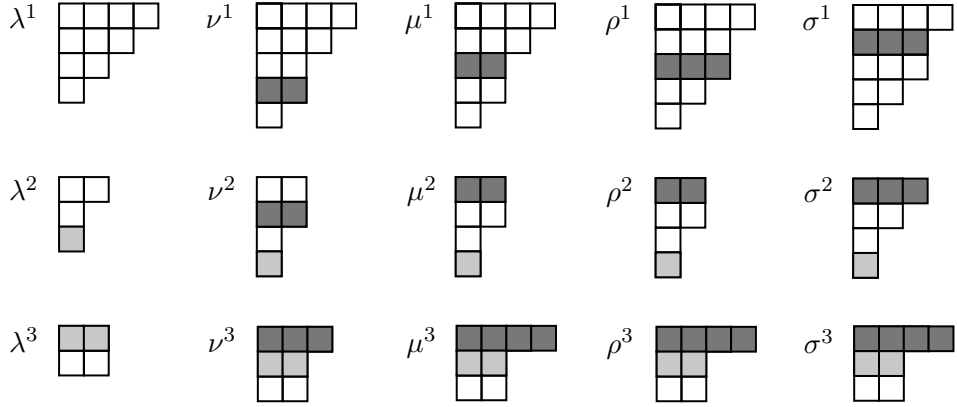,height=5.333cm}
\end{center}
\caption{Consider $\la^1$, $\la^2$, and $\la^3$ with $p_2 = 4$ and $p_3 = 4$.  
We insert 7 into $\la^1$, $\la^2$, and $\la^3$ to get $\nu^1$, $\nu^2$, and $\nu^3$, 
8 to get $\mu^1$, $\mu^2$, and $\mu^3$, 9 to get $\rho^1$, $\rho^2$, and $\rho^3$, and 
10 to get $\si^1$, $\si^2$, and $\si^3$. }
\label{figInsertMore2}
\end{figure}
\bigskip

Condition~(1) follows immediately from either case since we have only added $1$ to one part.  
Also note that we never add 1 to a row that already has length $p_i$ which implies condition (2).  

Finally, consider the selected parts of $\mu^1, \mu^2, ..., \mu^k$.  Let $i$ be as found above.  
For partitions  $\mu^{i+1}, ..., \mu^k$ we select the same part as in $\nu^{i+1}, ..., \nu^k$.  
In $\mu^i$ we select the part to which we added 1.  For partitions  $\mu^{1}, ..., \mu^{i-1}$ we 
select the part directly above the selected part of $\nu^{1}, ..., \nu^{i-1}$ but because of our 
choice of $i$ these selected parts are equal to the selected parts of $\nu^{1}, ..., \nu^{i-1}$.  
Therefore selected parts have the same length as those inserted and 
$A(\nu^1, \nu^2, ..., \nu^k) = A(\mu^1, \mu^2, ..., \mu^k) + 1$, implying condition (3).  
\end{proof}

\begin{proof}[Proof of Proposition~\ref{prop:insertion}] The two previous lemmas give the base
case and inductive step necessary to prove Proposition~\ref{prop:insertion} without the
uniqueness property.  All that is needed to complete the proof is to check the uniqueness of 
$\mu^1, \mu^2, ..., \mu^k$.  

Suppose $\mu^1, \mu^2, ..., \mu^k$ and $\nu^1, \nu^2, ..., \nu^k$ are both sequences satisfying
conditions (1), (2), and (3) of the proposition for some particular sequence $\la^1, \la^2, ...,
\la^k$ and integer $a \geq A(\la^1, \la^2, ..., \la^k)$.  

Since removing the selected parts of each sequence gives $\la^1, \la^2, ..., \la^k$, the sequences $\mu^1, \mu^2, ..., \mu^k$ and $\nu^1, \nu^2, ..., \nu^k$ must differ in a selected part.  Let $i$ be the largest index so that the selected part of $\mu^i$ and $\nu^i$ are not equal.  Since $i$ is the largest index where this happens, 
the selected parts of $\mu^i$ and $\nu^i$ must sit in the same row, say $j$.  Without loss of generality, ${\mu^i}_j > {\nu^i}_j$.

Our selection procedure now forces the selected part of $\mu^s$ to be greater than or equal to the selected part of $\nu^s$ for $s < i$, which gives us 
$$
A(\mu^1, \mu^2, ..., \mu^k) > A(\nu^1, \nu^2, ..., \nu^k).  
$$
However, both of these are equal to $a$ and so we have reached a contradiction. 
There cannot be a difference between the sequence of partitions $\mu^1, \mu^2, ..., \mu^k$ and 
the sequence of partitions $\nu^1, \nu^2, ..., \nu^k$.
\end{proof}

This proposition will be used repeatedly to establish the bijections in sections~\ref{secSym}.  
For convenience, we establish the following notation.  
Let
$$
\phi_{\{p_2, ..., p_k\}}(a; \la^1, \la^2, ..., \la^k) = 
\phi(a; \la^1, \la^2, ..., \la^k)= (\mu^1, \mu^2, ..., \mu^k)
$$
where $\mu^1, \mu^2, ..., \mu^k$ are the partitions uniquely defined by 
Proposition~\ref{prop:insertion}.  
Of course, $\phi$ is only defined for $\la^1, \la^2, ..., \la^k$ and $a$ such that 
$A(\la^1, \la^2, ..., \la^k) \leq a$.

The proof of Proposition~\ref{prop:insertion} 
gives us the following algorithm for insertion.
 
\begin{alg}
Let $\la^1, \la^2, ..., \la^k$ and $a$ be such that $A(\la^1, \la^2, ..., \la^k) \leq a$. 

First insert a part of the same length as the part selected from $\la^i$ when calculating $A(\la^1, \la^2, ..., \la^k)$ to $\la^i$ to obtain $\nu^i$.

Now we proceed recursively, adding one square at a time to $\nu^1, \nu^2, ..., \nu^k$ until we have inserted parts whose sum is $a$.  
To add one more box to the sequence of partitions:
\begin{itemize}
\item
If the selected part of $\nu^1$ is the first part, add one to this part.
\item
Otherwise,
find the partition $\nu^i$ with smallest index $i$ such that the selected part of $\nu^i$ is strictly less than the part above it or is strictly less than $p_i$ if it is the first part, and add one to this part.  
\end{itemize}
When we have added a total of $a$ boxes, let $\mu^1, \mu^2, ..., \mu^k$ be the resulting partitions.  We have
$$
\phi(a; \la^1, \la^2, ..., \la^k) = (\mu^1, \mu^2, ..., \mu^k) \, .
$$
\end{alg}

\subsection{Relationship between selection and insertion}
The last line of Proposition~\ref{prop:insertion} also shows that insertion  is reversible. Since the rows added by 
$\phi$ are those selected when calculating $A(\mu^1, \mu^2, ..., \mu^k)$ and 
$a = A(\mu^1, \mu^2, ..., \mu^k)$, we can undo insertion by using selection and removal.
It will be useful to formally note this consequence of Proposition~\ref{prop:insertion}.
\begin{cor}
\label{cor:phipsi}
Let $p_2$, $p_3$, ..., $p_k$ be integers.
\begin{enumerate}   
\item
For any sequence of $k$ partitions $\la^1, \la^2, ..., \la^k$ such that 
$$
f(\la^2) \leq p_2, \, f(\la^3) \leq p_3, \, ..., \, f(\la^k) \leq p_k 
$$
and integer $a$ such that $a \geq A(\la^1, \la^2, ..., \la^k; p_2, ..., p_k)$ we have
$$
\psi(\phi(a;\la^1, \la^2, ..., \la^k)) = (a;\la^1, \la^2, ..., \la^k).
$$
\item
For any sequence of $k$ partitions $\mu^1, \mu^2, ..., \mu^k$ such that 
$$
f(\mu^2) \leq p_2, \, f(\mu^3) \leq p_3, \, ..., \, f(\mu^k) \leq p_k 
$$
we have
$$
\phi(\psi(\mu^1, \mu^2, ..., \mu^k)) = (\mu^1, \mu^2, ..., \mu^k).
$$
\end{enumerate}
\end{cor}

\subsection{Iterative removal of selected parts}
\label{secIter}
\begin{figure}[hbt]
\begin{center}
\psfrag{1}{$\mu^1$}
\psfrag{2}{$\mu^2$}
\psfrag{3}{$\mu^3$}
\psfrag{4}{$\mu^4$}
\psfrag{5}{$\nu^1$}
\psfrag{6}{$\nu^2$}
\psfrag{7}{$\nu^3$}
\psfrag{8}{$\nu^4$}
\epsfig{file=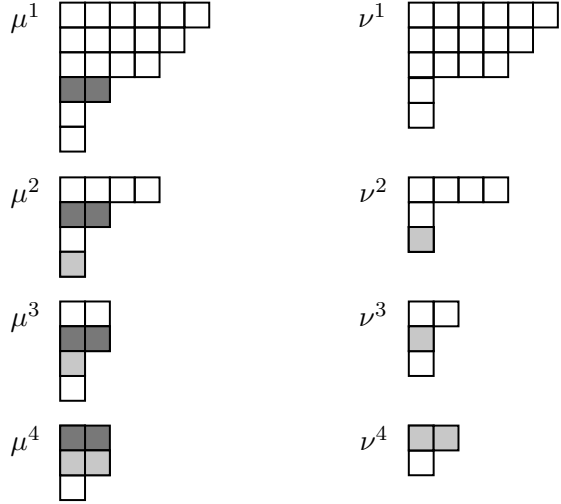,height=6.666cm}
\end{center}
\caption{For $p_2 = 4$, $p_3 = 2$, and $p_4 = 3$, we see that for $\psi_2(\mu^1, \mu^2, \mu^3, \mu^4) = (\nu^1, \nu^2, \nu^3, \nu^4)$.  Also note that $A(\mu^1, \mu^2, \mu^3, \mu^4) = 8 \geq A(\nu^1, \nu^2, \nu^3, \nu^4) = 4$.}
\label{figLemmaPsi}
\end{figure}

As a final remark, we note that if $\psi_2(\mu^1, \mu^2, ..., \mu^k) = (\la^1, \la^2, ..., \la^k)$, then
$$
f(\la^2) \leq p_2, \, f(\la^3) \leq p_3, \, ..., \, f(\la^k) \leq p_k \, .
$$
Therefore we can apply $\psi_1$ or $\psi_2$ to $\la^1, \la^2, ..., \la^k$ 
and in general we can reapply $\psi_2$, the removal of selected parts,  any number of times.  
The following lemma describes more precisely what happens to selected parts when $\psi_2$ is applied
more than once. 
See Figure~\ref{figLemmaPsi}.

\begin{lemma}
\label{lemma:selectbelow}
For any sequence of $k$ partitions $\mu^1, \mu^2, ..., \mu^k$ such that 
$$
f(\mu^2) \leq p_2, \, f(\mu^3) \leq p_3, \, ..., \, f(\mu^k) \leq p_k 
$$
the selected rows of $\psi_2(\mu^1, \mu^2, ..., \mu^k)$ are rows that sit 
{\emph strictly} below the selected rows of $\mu^1, \mu^2, ..., \mu^k$ in $\mu^1, \mu^2, ..., \mu^k$.

In particular we have 
$$
A(\psi_2(\mu^1, \mu^2, ..., \mu^k); p_2, p_3, ..., p_k) \leq A(\mu^1, \mu^2, ..., \mu^k; p_2, p_3, ..., p_k) \, .
$$
\end{lemma}

\begin{proof}
This follows by a simple inductive argument.  

Let $\psi_2(\mu^1, \mu^2, ..., \mu^k) = (\nu^1, \nu^2, ..., \nu^k)$.
In both $\mu^k$ and $\nu^k$ we select the first part.  However, the first part of $\nu^k$ is the second part of $\mu^k$ and so the result holds for $\mu^k$ and $\nu^k$. 

Moreover if the result is true for $\mu^i$ and $\nu^i$, and if we selected the $h$th row of $\mu^i$ and the $j$th row of $\nu^i$, then we have $h \leq j$.  This implies ${\mu^i}_h \geq {\nu^i}_j$.  Then the selected rows of $\mu^{i-1}$ and $\nu^{i-1}$ are $h + (p_i - {\mu^i}_h)$ and $j + (p_i - {\nu^i}_j)$ respectively and $h + (p_i - {\mu^i}_h) \leq j + (p_i - {\nu^i}_j)$ as desired.
\end{proof}

The procedures presented here are the main tools needed build the combinatorial proof
of~(\ref{eqn:genSchur}).  These procedures can also be used to obtain
other bijections as is shown in Appendix~A of~\cite{Bou:Thesis}.

%%%%%%%%%%%%%%%%%%%%%%%%%%%%%%%%%%%%%%%%%%%%%%%%%%%%%%%%%%%
%%%%                                                   %%%%
%%%%              Definition of (k,m)-rank             %%%%
%%%%                                                   %%%%
%%%%%%%%%%%%%%%%%%%%%%%%%%%%%%%%%%%%%%%%%%%%%%%%%%%%%%%%%%%

\section{Definition of $(k,m)$-rank}
\label{secDef}

In this section, we will define notion of $(k,m)$-rank for partitions with at least $k$ successive
$m$-Durfee rectangles.  (This generalizes the $(2,m)$-rank for partitions with at least two successive
$m$-Durfee rectangles found in~\cite{BP:RR}.)

\begin{figure}[hbt]
\begin{center}
\psfrag{1}{$\la^1$}
\psfrag{2}{$\la^2$}
\psfrag{3}{$\la^3$}
\psfrag{4}{$\al$}
\psfrag{a}{$\leftarrow N_1 \rightarrow$}
\psfrag{b}{$N_2$}
\psfrag{c}{$N_3$}
\psfrag{d}{$\, N_1$}
\epsfig{file=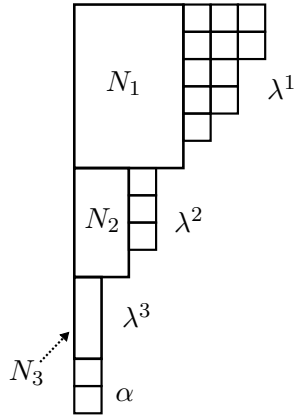,height=5.5cm}
\end{center}
\caption{Successive Durfee rectangles of width $N_1$, $N_2$, and $N_3$ and names for the partitions to the right, $\la^1$, $\la^2$, and $\la^3$, and below, $\al$, those Durfee rectangles.}
\label{figParName}
\end{figure}

First, given a partition $\la$ with $k$ successive $m$-Durfee rectangles, denote by $\la^i$ the
partition to the right of the $i$th $m$-Durfee rectangle and denote by $\al$ the partition below the
$k$th $m$-Durfee rectangle.  Moreover, let $N_1, N_2, ..., N_k$ denote the widths of the first $k$
successive $m$-Durfee rectangles.   Note that, for all $i$,  $\la^i$ has at most $N_i +m$ parts and, for
$i \geq 2$, the largest part of $\la^i$ is at most $N_{i-1} - N_i$.   See Figure~\ref{figParName}.  

Dyson's original definition of rank applies to a partition with (at least) one Durfee
square and compares the largest part of the $\la^1$ to the number of parts of $\al$.  Our
$(k,m)$-rank will compare parts to the right of the $k$ successive $m$-Durfee rectangles to the number
of parts of $\al$.  

Let $p_i = N_{i-1}-N_i$. Then, as we noted above, we have $f(\la^i) \leq p_i$ for $i \geq 2$. 
Therefore, we may apply selection to the sequence of partitions to the right of our Durfee rectangles,
$\la^1, \la^2, ..., \la^k$.  

\begin{figure}[hbt]
\begin{center}
\psfrag{s}{selected parts}
\psfrag{b}{}
\psfrag{c}{}
\epsfig{file=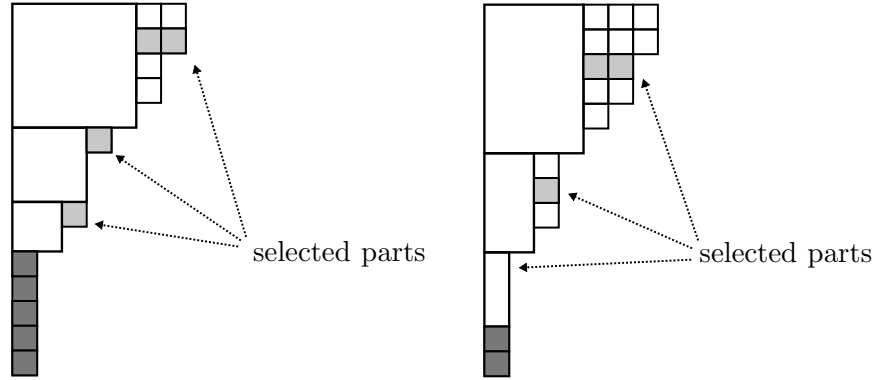,height=5cm}
\end{center}
\caption{For partition $\la = (7,7,6,6,5,4,3,3,3,2,1,1,1,1,1)$, we have $a_{3,0}(\la) = 2+1+1=4$, $b_{3,0}(\la) = 5$, and $r_{3,0}(\la) = 4-5 = -1$, while $a_{3,1}(\la) = 2+1+0 = 3$, $b_{3,1}(\la) = 2$, and $r_{3,1}(\la) = 3-2=1$.}
\label{figSelect}
\end{figure}

See Figure~\ref{figSelect} for two examples of this selection process.  On the left hand side, we
consider $\la$ with 3 successive Durfee squares.  First we select the first part of $\la^3$.  Next we
select the $1+(1-1) = 1$st part of $\la^2$ and we select the $1 + (2-1) =
2$nd part of $\la^1$.  On the right hand side, we consider $\la$ with 3 successive $1$-Durfee
rectangles.  First we select the first part of $\la^3$, then we select
the $1 + (1-0) = 2$nd part of $\la^2$ and finally the $2 + (2-1) = 3$rd part of $\la^1$.  

Notice that in these examples, the selected part of the partition $\la^i$ is never below the bottom row
of the $m$-Durfee rectangle sitting to its left.  This is true in general as stated by the following lemma.

\begin{lemma}
If the $j$th part of $\la^i$ has been selected, then $j  \leq N_i +m$.
\end{lemma}

\begin{proof}
We will prove the stronger statement that if the $j$th part of $\la^i$ has been selected then 
$j \leq 1 + N_i - N_k$.  

If $m \leq 0$, the $k$th successive $m$-Durfee square has non-zero height, and so its width is $N_k \geq 1 - m$, which gives us 
$$1 +N_i - N_k \leq N_i +m \, .$$
If $m > 0$, we get $N_k \geq 0$ so 
$$1 +N_i - N_k \leq N_i + 1 \leq N_i +m \, .$$  
Therefore the statement in the lemma follows from $j \leq 1 + N_i - N_k$.  

To show that $j \leq 1 + N_i - N_k$, we proceed by induction, starting with $\la^k$ and moving up to $\la^1$.  

We select the first row of $\la^k$ and have $1 = 1+ N_k -N_k$.  

If the $j$th row of $\la^i$ has been selected, we select the 
$j + p_i -{\la^i_j} = (j+ (N_{i-1} - N_i) - {\la^i}_j)$th
row of $\la^{i-1}$.  
Now our inductive hypothesis says that then $j  \leq 1 +N_i - N_k$.
Hence we see that 
$$
\aligned
j +  (N_{i-1} - N_i) - {\la^i}_j & \leq  1 + N_i - N_k + N_{i-1} - N_i - {\la^i}_j \\
 & \leq  1 + N_{i-1} - N_k
\endaligned
$$ as desired.
\end{proof}

Finally, we can give the definition of $(k,m)$-rank, $r_{k,m}(\la)$.  
\begin{defn}
For $k > 0$, consider a partition $\la$ with $k$ successive $m$-Durfee rectangles of width $N_1$, $N_2$,
..., $N_k$.  Let $p_i = N_{i-1}-N_i$. Also, let $\la^1, \la^2, ..., \la^k$ be the partitions to the right of
the Durfee rectangles and let $\al$ be the partition below the $k$th Durfee rectangle.  Define  
\begin{itemize}
\item
$a_{k,m}(\la) = A(\la^1, \la^2, ..., \la^k; p_2, p_3, .., p_k)$, the sum of the parts selected from $\la^1, \la^2, ..., \la^k$,
\item
$b_{k,m}(\la) = \ell(\al)$, the number of parts of $\al$ and
\item
$r_{k,m}(\la) = a_{k,m}(\la) - b_{k,m}(\la)$.
\end{itemize}
\end{defn}

In words, our definition of $(k,m)$-rank selects parts to the right of the $k$ successive $m$-Durfee rectangles of $\la$
and compares the total size of these parts to the number of parts below the Durfee rectangles.  In the case $k = 1$ and
$m=0$ this corresponds exactly to Dyson's original definition.  See Figure~\ref{figSelect} for examples.

%%%%%%%%%%%%%%%%%%%%%%%%%%%%%%%%%%%%%%%%%%%%%%%%%%%%%%%%%%%
%%%%                                                   %%%%
%%%%                    Symmetries                     %%%%
%%%%                                                   %%%%
%%%%%%%%%%%%%%%%%%%%%%%%%%%%%%%%%%%%%%%%%%%%%%%%%%%%%%%%%%%

\section{Symmetries}
\label{secSym}

Let $h(n,k,m,r)$ be the number of partitions of $n$ with $(k,m)$-rank equal to $r$.  
Similarly,  let $h(n,k,m,\leq r)$ be the number of partitions of $n$ with $(k,m)$-rank less than or equal to $r$ and 
let $h(n,k,m, \geq r)$ be the number of partitions of $n$ with $(k,m)$-rank greater than or equal to $r$.  

There are relationships between these numbers that generalize the symmetries used in Dyson's
proof of Euler's pentagonal number theorem.  These relationships will be proved in a completely
combinatorial way and, in the following section, they will be used to establish the generalized Schur
identities by simple algebraic manipulation.

Recall that $q_{k-1}(n)$ denotes the number of partitions with at most $k-1$ Durfee squares, so that
$p(n) - q_{k-1}(n)$ is the number of partitions with at least $k$ Durfee squares.  

The following two observations follow directly from our definitions since $(k,0)$-rank is only defined
for the set of partitions with $k$ non-empty Durfee squares, whereas when $m > 0$, $(k,m)$-rank is
defined for the set of all partitions.  

Unless otherwise stated we implicitly assume that $n,k,m,r \in \mathbb{Z}$ and $k > 0$.  

\begin{obs}[First observation]
\label{genObs1}
For $m = 0$,
$$
h(n,k,0, \leq r) + h(n,k,0,\geq r+1) = p(n) - q_{k-1}(n) \, .
$$
\end{obs}

\begin{obs}[Second observation]
\label{genObs2}
For $m > 0$,
$$
h(n,k,m, \leq r) + h(n,k,m,\geq r+1) = p(n) \, .
$$
\end{obs}

There are also more complicated relations between these numbers.  

\begin{theorem}[First symmetry]
\label{genSym1}
For any $r \in \zz$ and $m=0$,
$$
h(n,k,0,r) = h(n,k,0,-r) \, .
$$
\end{theorem}

\begin{theorem}[Second symmetry]
\label{genSym2}
For $m,r \in \zz$, if  $r > 0$ or if $m \leq 0$,
$$
h(n,k,m,\leq -r) = h(n -r- k(m+1)   ,k,m+2, \geq -r) \, .
$$
\end{theorem}

These symmetries generalize the symmetries in Dyson's proof and will be proved bijectively.  The first
will be proved by a bijection  generalizing conjugation and the second by a bijection generalizing
$d_r$, Dyson's map for his new symmetry.

\subsection{First symmetry}
\label{sec:Sym1}

To prove the first symmetry we are looking for a conjugation-like bijection that takes as its input a
partition with $(k,0)$-rank $r$ and outputs a partition with $(k,0)$-rank $-r$.  The following theorem
says that this bijection will not change the Durfee square structure of the partition and suggests that
we find a way to simply exchange the selected parts of $\la^1, ..., \la^k$, the partitions to the right
of the Durfee squares, with the first column of $\al$, the partition below the Durfee squares.  

\begin{theorem} 
\label{thm:sym1} 
For any integers $s, t \geq 0$, the number of partitions $\la$ of $n$ with $k$ successive Durfee
squares of widths $N_1$, $N_2$, ..., $N_k$ such that $a_{k,0}(\la) = s$ and $b_{k,0}(\la) = t$ is equal
to the number of partitions $\mu$ of $n$ with $k$ successive Durfee squares of widths $N_1$, $N_2$, ...,
$N_k$ such that $a_{k,0}(\mu) = t$ and $b_{k,0}(\mu) = s$.   
\end{theorem}

The following two corollaries follow immediately from this theorem.

\begin{cor}
For any $r \in \zz$, the number of partitions $\la$ of $n$ with $k$ successive Durfee squares of widths $N_1$, $N_2$, ..., $N_k$ such that $r_{k,0}(\la) = r$ is equal to the number of partitions $\mu$ of $n$ with $k$ successive Durfee squares of widths $N_1$, $N_2$, ..., $N_k$ such that $r_{k,0}(\mu) = -r$.
\end{cor}

\begin{cor}
\label{corSym1}
For any $r \in \zz$, the number of partitions $\la$ of $n$ with $k$ successive Durfee squares and $r_{k,0}(\la) = r$ is equal to the number of partitions $\mu$ of $n$ with $k$ successive Durfee squares and $r_{k,0}(\mu) = -r$.  
\end{cor}

Note that Corollary~\ref{corSym1} is exactly Theorem~\ref{genSym1}, the first symmetry.

To prove Theorem~\ref{thm:sym1} and its corollaries, we present a map,
$$
\mathfrak{C}^k: \mathcal{P} \smallsetminus \mathcal{Q}_{k-1} \rightarrow \mathcal{P} \smallsetminus \mathcal{Q}_{k-1} \, .
$$

Recall that 
$$
\mathcal{Q}_{k-1} \,  = \, \{\text{partitions with at most $k-1$ Durfee squares}\} \, ,
$$  
and so
$$
\mathcal{P} \smallsetminus \mathcal{Q}_{k-1} \,  = \, \{\text{partitions with at least $k$ Durfee squares}\} \, .
$$

\begin{alg}[Generalized Conjugation]
\label{algC}
Let $\la$ be a partition with at least $k$ Durfee squares.  

Let $\alpha$ be the partition below the $k$th Durfee square of $\la$ and $\la^1, \la^2, ..., \la^k$ be
the partitions to the right of these squares.  

Let $N_1, N_2, ..., N_k$ be the size of the $k$ successive Durfee squares and let 
$$
p_2 = N_1 - N_2, \, p_3 = N_2 - N_3, \, ..., \,  p_k = N_{k-1} - N_k .
$$

We iteratively remove selected parts  from $\la^1, \la^2, ..., \la^k$ by using $\psi_2$ $N_k$ times to obtain:
$$
{\psi_2}^{N_k}(\la^1, \la^2, ..., \la^k) = (\nu^1, \nu^2, ..., \nu^k) \, .
$$

As these rows are removed, we record the sum of the parts removed at each step.  We record this
information in a partition $\beta$ as follows. The partition $\beta$ is defined by giving its conjugate.  
$$
\aligned
\be^\pr_1       \,  = & \, \psi_1(\la^1, \la^2, ..., \la^k) \, , \\[3mm]
\be^\pr_2       \,  = & \, \psi_1(\psi_2(\la^1, \la^2, ..., \la^k)) \, , \\[3mm]
                      & \vdots \\[2mm] 
\be^\pr_{N_k-1} \,  = & \, \psi_1({\psi_2}^{N_k-2}(\la^1, \la^2, ..., \la^k)) \, , \\[3mm]
\be^\pr_{N_k}   \,  = & \, \psi_1({\psi_2}^{N_k-1}(\la^1, \la^2, ..., \la^k)) \, .
\endaligned
$$

Next we insert $\alpha^\pr_{N_k}, \alpha^\pr_{N_k-1}, ..., \alpha^\pr_1$ into $\nu^1, \nu^2, ..., \nu^k$ in that order giving us
$$
\phi(\al^\pr_1; \phi(\al^\pr_2; ... \phi(\al^\pr_{N_k-1}; \phi(\al^\pr_{N_k}; \nu^1, ..., \nu^k)) ...))  = (\mu^1, ..., \mu^k) \, .
$$

Let $\mu$ be a new partition defined by having
\begin{itemize} 
\item
$k$ successive Durfee squares of widths $N_1 , N_2 , ..., N_k $, 
\item
$\mu^1, \mu^2..., \mu^k$ to the right of these squares, and
\item
$\be$ below the $k$th Durfee square.  
\end{itemize}
Then $\mathfrak{C}^k(\la) = \mu$.
\end{alg}

We will also write $\mathfrak{C}$ for $\mathfrak{C}^k$.  

The map $\mathfrak{C}$ consists of using the tools from Section~\ref{secTools} to first remove selected
parts from the partitions to the right of the Durfee squares $N_k$ times and record the sum of the
removed parts as the $N_k$ columns of the partition $\beta$ (to be place below the Durfee squares). 
Second, we insert the $N_k$ columns of the partition $\al$ (that was below the Durfee squares) into the
partitions to the right of the Durfee squares (from which we just removed parts).  As such we are
exchanging some parts to the right of the Durfee squares with the parts below the Durfee squares. 

\begin{figure}[hbt]
\begin{center}
\psfrag{1}{$\la$}
\psfrag{2}{$\mathfrak{C}^2(\la)$}
\psfrag{a}{$\nu^1$}
\psfrag{b}{$\nu^2$}
\psfrag{c}{$\nu^3$}
\psfrag{d}{$\nu^4$}
\psfrag{e}{$\al^\pr$}
\psfrag{f}{$\be^\pr$}
\epsfig{file=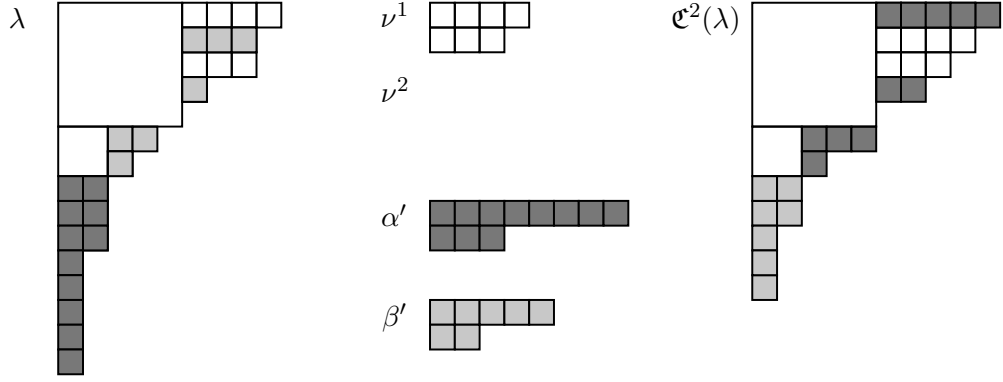,height=5cm}
\end{center}
\caption{Applying $\mathfrak{C}^2$ to the partition $\la = (9,8,8,6,5,4,3,2,2,2,1,1,1,1,1)$
gives $\mathfrak{C}^2(\la) = (10,9,8,7,5,5,3,2,2,1,1,1)$.  Intermediate steps are $\nu^1, \nu^2, \al^\pr$, and $\be^\pr$ as shown.}
\label{figCongEg2}
\end{figure}

\begin{figure}[hbt]
\begin{center}
\psfrag{1}{$\la$}
\psfrag{2}{$\mathfrak{C}^4(\la)$}
\psfrag{a}{$\nu^1$}
\psfrag{b}{$\nu^2$}
\psfrag{c}{$\nu^3$}
\psfrag{d}{$\nu^4$}
\psfrag{e}{$\al^\pr$}
\psfrag{f}{$\be^\pr$}
\epsfig{file=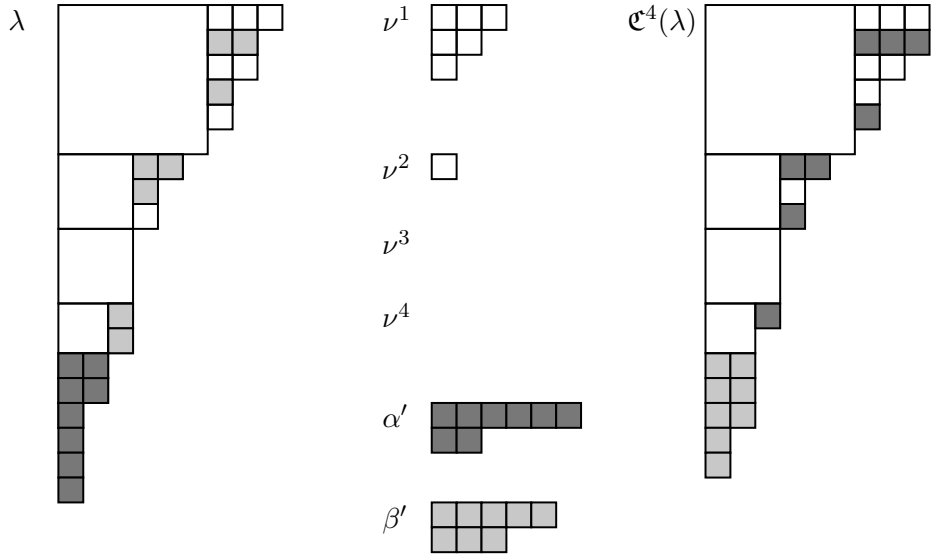,height=7.3333cm}
\end{center}
\caption{Applying $\mathfrak{C}^4$ to $\la = (9,8,8,7,7,6,5,4,4,3,3,3,3,3,2,2,1,1,1,1)$
gives $\mathfrak{C}^4(\la) = (9,9,8,7,7,6,5,4,4,3,3,3,3,2,2,2,2,1,1)$.  Intermediate steps are $\nu^1, \nu^2, \nu^3, \nu^4, \al^\pr$, and $\be^\pr$ as shown.}
\label{figCongEg1}
\end{figure}

Because of the way the maps $\phi$ and $\psi$ are defined, this exchange of parts is well-defined and is
an involution.   
Before proving this we will give two examples of applications of $\mathfrak{C}$.   These are found in figures~\ref{figCongEg2} and~\ref{figCongEg1}.  

\begin{proof}[Proof of Theorem~\ref{thm:sym1}]

Let $\la$ be a partition with 
$$
\aligned
a_{k,m}(\la) & =  s \, \text{ and}\\
b_{k,m}(\la) & =  t \, .
\endaligned
$$

To prove Theorem~\ref{thm:sym1}, we will show that $\mathfrak{C}$ is an involution that exchanges $a_{k,m}$ and $b_{k,m}$ while
preserving the Durfee square structure of $\la$.  

First we want to show that $\mathfrak{C}$ is well-defined.  

For each $2 \leq i \leq k$, $\la^i$ fits to the right of the $i$th Durfee rectangle and below the $(i-1)$st Durfee rectangle.  As a consequence, its largest part satisfies $f(\la^i) \leq N_i - N_{i-1} = p_i$.  

Therefore, we may select parts from and apply the maps $\psi_1$ and $\psi_2$ to $\la^1, \la^2, ..., \la^k$.  
Moreover, as we remarked in of Section~\ref{secIter}, the iterated applications we do here are also fine.  

Next we want to check that $\be$ is a partition and that $\be$ fits belo the $k$th Durfee square.  
Lemma~\ref{lemma:selectbelow} tells us that
$$
\psi_1(\psi_2(\mu^1, \mu^2, ..., \mu^k)) \leq \psi_1(\mu^1, \mu^2, ..., \mu^k) 
$$
for any $\mu^1, \mu^2, ..., \mu^k$.
Therefore, $$
\psi_1(\psi_2^j(\la^1, \la^2, ..., \la^k)) \leq \psi_1(\psi_2^{j-1}(\la^1, \la^2, ..., \la^k)) 
$$ for all $j \geq 1$.
In other words, $\be^\pr_{j+1} \leq \be^\pr_j$.  
Since $\be^\pr_1 \geq \be^\pr_2 \geq ... \geq \be^\pr_{N_k}$, we see that $\be$ a partition.  
Also $\be^\pr$ has at most $N_k$ parts which implies that $f(\be) \leq N_k$.  
Therefore, we can place $\be$ below the $k$th Durfee square whose size is $N_k$.   

We must also show that $\al^\pr$ can be inserted into the partitions to the right of the Durfee squares.  
Consider 
$$
{\psi_2}^{N_k}(\la^1, \la^2, ..., \la^k) = (\nu^1, \nu^2, ..., \nu^k) \, .
$$
Since we have simply removed parts, we have 
$$
f(\nu^2) \leq p_2, \, f(\nu^3) \leq p_3, \, ..., \, f(\nu^k) \leq p_k .
$$

Moreover, for $1 \leq i \leq k$, $\nu^i$ has at most $N_i - N_k$ parts.   (In particular, $\nu^k$ is empty.)  
This means that when we select parts from $\nu^1, \nu^2, ..., \nu^k$, we select the $(1 + N_i - N_k)$th part of $\nu^i$ which is always empty.  As a consequence, $A(\nu^1, \nu^2, ..., \nu^k; p_2, p_3, ..., p_k) = 0$.

Therefore we can insert ${\al^\pr}_{N_k} \geq 0$ into $\nu^1, \nu^2, .., \nu^k$. 
As well, since 
$
A(\phi({\al^\pr}_i; ...)) = {\al^\pr}_i \, 
$
and
$
{\al^\pr}_i \leq {\al^\pr}_{i-1}
$
we can insert ${\al^\pr}_{i-1}$ into 
$\phi({\al^\pr}_i; ... \phi(\al^\pr_{N_{k-1}}; \phi(\al^\pr_{N_k}; \nu^1, ..., \nu^k)) ...)$.  

Finally, each of these insertions adds at most one part to each partition and does not give partitions whose largest parts are greater than $p_i$.  Therefore
$$
f(\mu^2) \leq p_2, \, f(\mu^3) \leq p_3, \, ..., \, f(\mu^k) \leq p_k
$$
and, for $1 \leq i \leq k$,  $\mu^i$ has at most $(N_i -N_k) + N_k = N_i$ parts.  
Each $\mu^1, \mu^2, ..., \mu^k$ can be inserted to the right of each of the first $k$ Durfee 
squares and this shows that $\mathfrak{C}$ is well-defined.

\bigskip

To see that $\mathfrak{C}$ is an involution, we simply use the relationship between $\phi$ and $\psi$ summarized 
in Corollary~\ref{cor:phipsi}.  Say $\mathfrak{C}(\la) = \mu$ with $\al$, $\be$, $(\la^1, \la^2, ..., \la^k)$, 
and $(\mu^1, \mu^2, ..., \mu^k)$ as in the definition of Procedure~\ref{algC}.  
We will apply $\mathfrak{C}$ to $\mu$.  

Applying $\psi$ undoes the insertions done by $\phi$, and 
we get
$$
\aligned
\psi_1(\mu^1, \mu^2, ..., \mu^k) & = \psi_1(\phi(\al^\pr_1; ...)) = \al^\pr_1 \, \text{ and}\\
\psi_1(\psi_2(\mu^1, \mu^2, ..., \mu^k)) & = \psi_1(\psi_2(\phi(\al^\pr_1; \phi(\al^\pr_2; ...)))) = \al^\pr_2 \, .
\endaligned
$$
Similarly,
$$
\aligned
\psi_1({\psi_2}^2(\mu^1, \mu^2, ..., \mu^k)) & =  \al^\pr_3 \, , \\
                   &    \vdots \\ 
\psi_1({\psi_2}^{N_k-2}(\mu^1, \mu^2, ..., \mu^k)) & =  \al^\pr_{N_k-1} \, ,\\
\psi_1({\psi_2}^{N_k-1}(\mu^1, \mu^2, ..., \mu^k)) & =  \al^\pr_{N_k} \, ,
\endaligned
$$
and
$$
{\psi_2}^{N_k}(\mu^1, \mu^2, ..., \mu^k))  = (\nu^1, \nu^2, .., \nu^k) \, .
$$

Next we insert $\be^\pr_{N_k}, \be^\pr_{N_k-1}, ..., \be^\pr_1$ into $\nu^1, \nu^2, .., \nu^k$.  
Since these are the parts originally removed by $\psi$ from $\la^1, \la^2, ..., \la^k$ to give 
$\nu^1, \nu^2, .., \nu^k$, we have
$$
\phi(\be^\pr_1; \phi(\be^\pr_2; ... \phi(\be^\pr_{N_k-1}; \phi(\be^\pr_{N_k}; \nu^1, ..., \nu^k)) ...))  
= (\la^1, \la^2, ..., \la^k) \, .
$$
Since $\al$ goes below the $k$th Durfee square and $\la^1$, $\la^2$, ..., $\la^k$ to the right of the Durfee 
squares, we get $\mathfrak{C}(\mu) = \la$.  
A second application of  
$\mathfrak{C}$ undoes the first application
and indeed $\mathfrak{C}$ is an involution.

Finally, we note that 
$$
a_{k,m}(\mu) = A(\mu^1, \mu^2, ..., \mu^k; p_2, p_3, ..., p_k) = \psi_1(\mu^1, \mu^2, ..., \mu^k) = \al^\pr_1 = b_{k,m}(\la) = t
$$
and since $\mathfrak{C}$ is an involution
$b_{k,m}(\mu) = a_{k,m}(\la) = s$.
\end{proof}

\subsection{Second Symmetry}

The second symmetry will follow from the following theorem.  It will be proved by giving a bijection that uses
insertion to generalize Dyson's map.  

\begin{theorem}
\label{thm:sym2}
For any $m,r,t \in \zz$ such that $t \geq 0$, 
the number of partitions $\la$ of $n$ with $k$ successive $m$-Durfee 
rectangles, of {\em non-zero} widths $N_1, N_2, ..., N_k$, with $b_{k,m}(\la) = t$ 
and $r_{k,m}(\la) \leq -r$ is equal to the number of partitions $\mu$ of $n-r-k(m+1)$ 
with $k$ successive $(m+2)$-Durfee rectangles, of widths $N_1-1$, $N_2-1$, ..., $N_k-1$, 
with $a_{k,m+2}(\mu) = t-r$ and $b_{k,m+2}(\mu) \leq t$.  
\end{theorem}

Before proving our theorem, note that the following two corollaries follow immediately.  

\begin{cor}
\label{cor:sym2}
For any $m,r \in \zz$, the number of partitions $\la$ of $n$ with $k$ successive
$m$-Durfee rectangles, of {\em non-zero} widths $N_1$, $N_2$, ..., $N_k$, with $r_{k,m}(\la) \leq -r$ is equal to
the number of partitions $\mu$ of $n-r-k(m+1)$ with~$k$ successive $(m+2)$-Durfee rectangles, of widths $N_1-1$,
$N_2-1$, ..., $N_k-1$, with~$r_{k,m+2}(\mu) \geq -r$.
\end{cor}

In Theorem~\ref{thm:sym2} and Corollary~\ref{cor:sym2}, one side of the identity requires {\em non-zero} width Durfee
rectangles while on the other side zero width Durfee rectangles are allowed. 
(Durfee rectangles of height zero are never allowed as stated in the definition of Durfee rectangles.)   
There are two situations in which the widths of the $k$ successive $m$-Durfee rectangles 
are forced to be non-zero.  

First, since we require $m$-Durfee rectangles to have non-zero height, the width of the
rectangles is at least $1-m$.  When $m \leq 0$, $1 -m > 0$ and so the width if forced to be non-zero.    

Second, if $r_{k,m}(\la) \leq -r$, then we must have $b_{k,m}(\la) \geq r$.  
If in addition $r > 0$, then $b_{k,m}(\la) > 0$.  Since $b_{k,m}(\la)$
is the size of the first column of $\al$ the partition which sits below the $k$th successive $m$-Durfee rectangle,
notice that this $m$-Durfee rectangle must have non-zero width.  Therefore if $r_{k,m}(\la) \leq -r$ and
$r > 0$, all $m$-Durfee rectangles of $\la$ have non-zero width.  

These cases give the following corollary which is Theorem~\ref{genSym2}.  

\begin{cor}
\label{corSym2}
For any $m,r \in \zz$ such that $m \leq 0$ or $r > 0$, the number of partitions
$\la$ of $n$ with $r_{k,m}(\la) \leq -r$ is equal to the number of 
partitions $\mu$ of $n-r-k(m+1)$ with $r_{k,m+2}(\mu)
\geq -r$.
\end{cor}

Note that Corollary~\ref{corSym2} is Theorem \ref{genSym2}, the second symmetry.

To prove Theorem~\ref{thm:sym2} and its corollaries, we present a family of maps,
$$
\mathfrak{D}^{k,m}_r: \mathcal{A} \rightarrow \mathcal{B} \, ,
$$
between the following two sets:
$$\aligned
\mathcal{A} \, = & \, \{\text{partitions with $k$ successive $m$-Durfee rectangles of {\em non-zero}} \\
                 & \ \ \ \text{width with $(k,m)$-rank at most $-r$} \} \, ,\\
\mathcal{B} \, = & \, \{\text{partitions with $(k,m+2)$-rank at least $-r$}\} \, .
\endaligned
$$

\begin{alg}[Generalized Dyson's map]
Let $\la$ be a partition with $r_{k,m}(\la) \leq -r$.  

Let $\alpha$ be the partition below the $k$th successive $m$-Durfee rectangle and $\la^1, \la^2, ..., \la^k$ be the partitions to the right of the rectangles.

Let $N_1, N_2, ..., N_k$ be the widths of the $k$ successive $m$-Durfee rectangles and let 
$$
p_2 = N_1 - N_2, \, p_3 = N_2 - N_3, \, ..., \,  p_k = N_{k-1} - N_k .
$$

Say $\ell(\alpha) = t$.  Then we obtain $k$ new partitions $\mu^1, \mu^2, ..., \mu^k$ by applying the insertion lemma so that
$$
\phi(t-r; \la^1, ..., \la^k) = (\mu^1, ..., \mu^k) \, .
$$

Remove the first column from $\alpha$ (or equivalently subtract 1 from each part) to get a partition $\beta$.

Let $\mu$ be a new partition defined by having
\begin{itemize} 
\item
$k$ successive $(m+2)$-Durfee rectangles of widths $N_1 - 1, N_2 - 1, ..., N_k -1$, 
\item
$\mu^1, \mu^2..., \mu^k$ to the right of these rectangles, and
\item
$\be$ below the $k$th rectangle.  
\end{itemize}

Then $\mathfrak{D}^{k,m}_r(\la) = \mu$.
\end{alg}

When $k$ and $m$ are clear from context we will write $\mathfrak{D}_r$.  

The essence of Dyson's map is to the remove the first column of a partition and, after adding or removing some boxes,
make it the first row of the partition.  Our map $\mathfrak{D}_r$, removes the first column of the partition below
the successive Durfee rectangle and inserts it (minus $r$ boxes) into the sequence of partitions to the right of the Durfee
rectangles.  To do this, the shape of the Durfee rectangle is modified to be one row taller and one column narrower.

\begin{figure}[hbt]
\begin{center}
\psfrag{1}{$\la$}
\psfrag{2}{$\mathfrak{D}_0(\la)$}
\epsfig{file=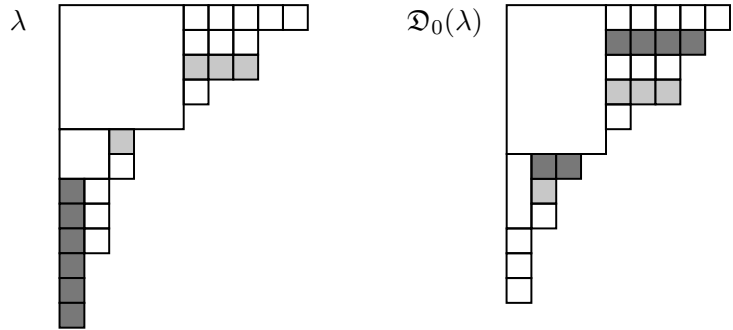,height=4.333cm}
\end{center}
\caption{Applying $\mathfrak{D}_0^{2,0}$ to $\la = (10,8,8,6,5,3,3,2,2,2,1,1,1)$.}
\label{figDysonEg1}
\end{figure}

\begin{figure}[hbt]
\begin{center}
\psfrag{1}{$\la$}
\psfrag{2}{$\mathfrak{D}_1(\la)$}
\epsfig{file=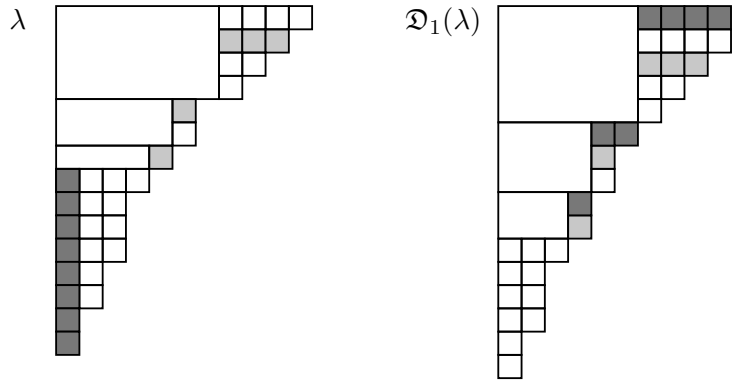,height=5cm}
\end{center}
\caption{Applying $\mathfrak{D}_1^{3,-3}$ to $\la = (11,10,9,8,6,6,5,4,3,3,3,2,2,1,1)$.}
\label{figDysonEg2}
\end{figure}

\begin{figure}[hbt]
\begin{center}
\psfrag{1}{$\la$}
\psfrag{2}{$\mathfrak{D}_3(\la)$}
\epsfig{file=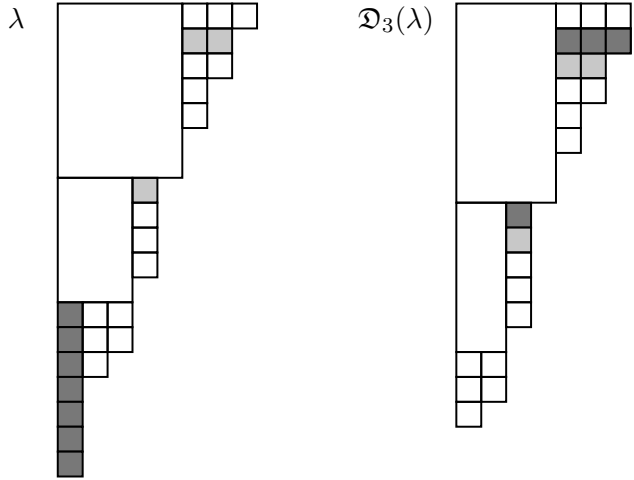,height=6.333cm}
\end{center}
\caption{Applying $\mathfrak{D}_3^{2,2}$ to $\la = (8,7,7,6,6,5,5,4,4,4,4,3,3,3,2,1,1,1,1)$.}
\label{figDysonEg3}
\end{figure}

We will give three examples of applications of $\mathfrak{D}_r$ before giving the proof that
$\mathfrak{D}_r$ is well-defined and gives a bijection between $\mathcal{A}$ and $\mathcal{B}$ that has the desired properties.  See Figures~\ref{figDysonEg1},~\ref{figDysonEg2}, and~\ref{figDysonEg3}.

\bigskip

\bigskip

\bigskip

\begin{proof}[Proof of Theorem~\ref{thm:sym2}]
To prove theorem~\ref{thm:sym2}, we show that $\mathfrak{D}^{k,m}_r$ is a bijection that changes rank and the other
statistics of the partitions appropriately.  
 
Consider a partition $\la$ with $k$ successive $m$-Durfee rectangles of widths $N_1$, $N_2$, ..., $N_k$ such that 
$r_{k,m}(\la) \leq -r$.  
Say 
$$
\aligned
a_{k,m}(\la) & =  s \, , \\
b_{k,m}(\la) & = t \, .
\endaligned
$$

First, we note that $\mu$ may have $(m+2)$-Durfee rectangles of width $N_1 -1, N_2 -1, ..., N_k-1$ since none 
of these integers are negative.  (If $m=-2$, none of these can be zero either since the $-2$-Durfee rectangles  of $\la$
must have nonzero height.

Next we want to apply our insertion procedure and so we must verify that the conditions of
Proposition~\ref{prop:insertion} are satisfied.   For each $2 \leq i \leq k$, $\la^i$ fits to the right of the
$i$th $m$-Durfee rectangle and below the $(i-1)$st $m$-Durfee rectangle.  As a consequence, its largest part satisfies
$f(\la^i) \leq N_i - N_{i-1} = p_i$ for $i \geq 2$ and $\ell(\la^i) \leq N_i +m$ for all $i \geq 1$.  

We want to insert $t - r$ into $\la^1, \la^2, ..., \la^k$.  Since $r_{k,m}(\la) = s-t \leq -r$, we get
$A(\la^i, \la^2, ..., \la^k; p_2, p_3, ..., p_k) = a_{k,m}(\la) = s \leq t-r$ which is the required condition.  

Now applying the lemma gives partitions $\mu^1, \mu^2, ..., \mu^k$ by inserting one part into 
each of $\la^1, \la^2,
..., \la^k$.  Since we doing so gives $f(\mu^i) \leq p_i = (N_i-1) - (N_{i-1}-1)$  and 
$\ell(\mu^i) \leq N_i +m +1 =
(N_i-1)+(m+2)$.  Hence  we see that $\mu^1, \mu^2, ..., \mu^k$
will fit to the right of the first $k$ successive $(m+2)$-Durfee squares of width $N_1 -1, N_2 -1, ..., N_k-1$.  

Finally, the largest part of $\be$ is one less that the largest part of $\al$ and so $\be$ fits under the
$k$th successive $(m+2)$-Durfee rectangle of $\mu$.  (If $\al$ was empty, $\be$ will be empty as well and will still
fit.) 

We may conclude that $\mu$ is a well-defined partition.  

We must show that $\mathfrak{D}_r$ is a bijection.  Notice that $\mathfrak{D}_r$ is reversible since
$$
\aligned
\psi_1(\mu^1, \mu^2, ..., \mu^k) & = t-r \, , \text{and}\\
\psi_2(\mu^1, \mu^2, ..., \mu^k) & = (\la^1, \la^2, ..., \la^k) \, .
\endaligned
$$
Hence we know $\la^1$, $\la^2$, ..., $\la^k$ and from this we can also recover 
$\al$ since we know $\be$ and since $\ell(\al) = t = \psi_1 +r$. 

To show that $\mathfrak{D}_r$ is surjective for any partition $\mu$ with $k$ successive $(m+2)$-Durfee 
rectangles of widths $N_1 -1, N_2-1, ..., N_k-1$ with $r_{k,m+2}(\mu) \geq -r$
we will construct a partition $\la$ such that
$\mathfrak{D}_r(\la) = \mu$.  
Since $r_{k,m+2}(\mu) \geq -r$, there is some $t\geq 0$ such that, $\mu$ has $a_{k,m+2}(\mu) = t-r$ and $b_{k,m+2}(\mu) \leq t$.  
Let $\mu^1, \mu^2, ..., \mu^k$ be the partitions to the right of the $k$ successive $(m+2)$-Durfee rectangles and let $\be$ be the partition below the $k$th successive $(m+2)$-Durfee rectangle.  

If we apply $\psi$ as above, we obtain partitions  $\la^1, \la^2, ..., \la^k$ of the appropriate size to put to the right of $m$-Durfee rectangles of width $N_1, N_2, ..., N_k$.
We can put a column to height $t = a_{k,m+2}(\mu) +r$ in front of $\be$ since $b_{k,m+2}(\mu) \leq t$
This partition $\al$ fits below the $k$th $m$-Durfee rectangle of width $N_k$.  
This gives a partition $\la$.

If we apply $\mathfrak{D}_r$ to $\la$ we are simply reversing the steps described above and so we get
$\mathfrak{D}_r(\la) = \mu$. This shows that $\mathfrak{D}_r$ is surjective onto the set of partitions with $k$
successive $(m+2)$-Durfee rectangles.   

We conclude that $\mathfrak{D}_r$ is indeed a bijection 

It remains to check the required properties of $\mu$.  
We note that:
\begin{itemize}
\item
by definition $\mu$ is a partition with $k$ successive $(m+2)$-Durfee rectangles of widths \\ $N_1 -1, N_2 -1, ..., N_k-1$,
\item
since we inserted $t-r$ into $\la^1, \la^2, ..., \la^k$ we get
$$
a_{k,m}(\mu) = A(\mu^1, \mu^2, ..., \mu^k; p_2, p_3, ..., p_k) =t-r \, ,
$$
\item
since $\ell(\al) = t$ we get
$$
b_{k,m}(\mu) = \ell(\be) \leq \ell(\al) = t \, , and
$$  
\item
if $\la$ is a partition of $n$, $\mu$ is a partition of $n -r -k(m+1)$ 
since we remove a column of height $t$ from $\al$ to get $\be$,  
insert $t-r$ into $\la^1, \la^2, ..., \la^k$ to get $\mu^1, \mu^2, ..., \mu^k$, and 
reduce the number of boxes in each of the $k$ successive Durfee rectangle by $m+1$.  
\end{itemize}
For the two corollaries we note that $r_{k,m+2}(\mu) = a_{k,m+2}(\mu) - b_{k,m+2}(\mu) \geq t-r -t = -r$.
\end{proof}

%%%%%%%%%%%%%%%%%%%%%%%%%%%%%%%%%%%%%%%%%%%%%%%%%%%%%%%%%%%
%%%%                                                   %%%%
%%%%   Algebraic derivation of the generalized Schur   %%%%
%%%%                                                   %%%%
%%%%%%%%%%%%%%%%%%%%%%%%%%%%%%%%%%%%%%%%%%%%%%%%%%%%%%%%%%%

\section{Algebraic derivation of the generalized Schur identity}
\label{secAlg}

We can now complete the proof of the generalized Schur identities.  
We proceed in a fashion similar to the algebraic steps of Dyson's proof of 
Euler's pentagonal number theorem using
the two observations and the two symmetries of the previous section.  

For every $j \in \mathbb{Z}$ let
$$\aligned
a_j \, & = \, h\left(n- jr - \frac{j(j-1)}{2} - k(jm+j^2),k, m+2j, \leq -r -j\right),\\[3mm]
b_j \, & = \, h\left(n- jr - \frac{j(j-1)}{2} - k(jm+j^2),k, m+2j, \geq -r -j+1\right).
\endaligned
$$

In this notation, for $m +2j > 0$, the second observation,~\ref{genObs2}, gives us
$$
a_j + b_j = p(n -jr -\frac{j(j-1)}{2} - k(jm+j^2)) \, .
$$  
For either $r +j > 0$ or for $m+2j \leq 0$, the second symmetry,~\ref{genSym2}, gives us
$$
a_j = b_{j+1} \, .
$$ 

Applying these multiple times we get
$$ \aligned
& h(n,k,m,\leq -r) \, = \, a_0 \,  = \, b_1 \\[3mm]
& \qquad \, =  \, b_1 + (a_1 - b_2) - (a_2 - b_3) + (a_3 - b_4) - \ldots \\[3mm]
& \qquad \, =  \, (b_1 + a_1) - (b_2 + a_2) + (b_3 + a_3)  - (b_4 + a_4) + \ldots \\[2mm]
& \qquad \, = \, \sum_{j=1}^\infty (-1)^{j-1} p(n -jr -\frac{j(j-1)}{2} - k(jm+j^2)) \, .
\endaligned
$$
This  last identity holds if either $m \geq 0$ and $r \geq 1$ or if $m = r = 0$.  

In terms of the generating functions
$$
\aligned
H_{k,m,\leq r} (q) \, &:= \, \sum_{n=0}^\infty \, h(n,k,m,\leq r) \, q^n \,, \text{ and}\\
H_{k,m,\geq r} (q) \, &:= \, \sum_{n=0}^\infty \, h(n,k,m,\geq r) \, q^n \,,
\endaligned
$$
this gives, if either $m \geq 0$ and $r \geq 1$ or if $m = r = 0$, 
\begin{equation*}
\label{eqn:genfuncrank}
H_{k, m,\leq -r} (q) \, = \,
\frac{1}{(q)_\infty} \ \sum_{j=1}^\infty (-1)^{j-1}
\, q^{jr +\frac{j(j-1)}{2} + k(jm+j^2)}\,.
\end{equation*}
In particular, we have:
$$
\aligned
H_{k, 0,\leq 0} (q) \, & = \, \frac{1}{(q)_\infty} \ \sum_{j=1}^\infty (-1)^{j-1}
\, q^{\frac{j(j-1)}{2} + kj^2}\,, \text{ and} \\
H_{k, 0,\le -1} (q) \, & = \, \frac{1}{(q)_\infty} \ \sum_{j=1}^\infty (-1)^{j-1}
\, q^{\frac{j(j+1)}{2} + kj^2}\,.
\endaligned
$$
From the first symmetry~\ref{genSym1} and the first observation~\ref{genObs1} we note that
$$H_{k,0,\leq 0} (q) + H_{k, 0,\leq -1}(q) \, = \,
H_{k, 0,\leq 0} (q) + H_{k, 0,\geq 1}(q)$$
is the generating function for partitions with at least $k$ successive Durfee squares.  
We conclude:
$$\aligned &
\frac{1}{(q)_\infty} \ \sum_{j=1}^\infty (-1)^{j-1}
\, q^{\frac{j(j-1)}{2} + kj^2} \, + \,
\frac{1}{(q)_\infty} \ \sum_{j=1}^\infty (-1)^{j-1}
\, q^{\frac{j(j+1)}{2} + kj^2}\, \\
& \qquad = \,
\frac{1}{(q)_\infty} \ -  \sum_{n_1 = 0}^\infty \cdots \sum_{n_{k-1} = 0}^\infty \frac{q^{N_1^2+N_2^2+\dots+N_{k-1}^2}}{(q)_{n_1}(q)_{n_2}\dots(q)_{n_{k-1}}}
\endaligned
$$
which implies the generalized Schur identities~(\ref{eqn:genSchur}) and completes our
proof of the generalized \RR identities~(\ref{eqn:genRR}).

%%%%%%%%%%%%%%%%%%%%%%%%%%%%%%%%%%%%%%%%%%%%%%%%%%%%%%%%%%%
%%%%                                                   %%%%
%%%%                    Connections                    %%%%
%%%%                                                   %%%%
%%%%%%%%%%%%%%%%%%%%%%%%%%%%%%%%%%%%%%%%%%%%%%%%%%%%%%%%%%%

\section{Connections to other work}
\label{secConn}

\subsection{Dyson's rank and proof of Euler's pentagonal number theorem}

As mentioned in the introduction, this proof of the 
\RR identities follows the general form of Dyson's proof of
Euler's pentagonal number theorem.   More specifically, our proof of the generalized Schur identities~(\ref{eqn:genSchur}) is a Dyson-style proof with a modified rank.   

We generalized Dyson's rank by
defining $(k,m)$-rank; his rank is our $(1,0)$-rank.   The algebraic steps used to deduce the generalized Schur
identities are the same as those used to deduce Euler's pentagonal number theorem.  Moreover, our symmetries and
corresponding bijections, $\mathfrak{C}^k$ and $\mathfrak{D}^{k,m}_r$  generalize conjugation and Dyson's map.  
More precisely, in the case $k =1$, we have:
$$
\aligned
r_{1,m}(\la) & = r(\la) -m \, , \\[2mm]
\mathfrak{C}^1 & = \text{usual conjugation}, \text{ and}\\[2mm]
\mathfrak{D}^{1,m}_{r} & = d_{-r-m} \, .
\endaligned
$$

This is not the first generalization of Dyson's rank that has been used to prove the \RR identities.  
The notion of
successive rank can also be used to give a combinatorial proof of the
Rogers-Ramanujan identities and their generalizations by a sieve argument
(see~\cite{Andr:sieves,Andr:hook,Bres:sieve}).   
However, this proof does not use the notion of successive Durfee squares but rather involves a different combinatorial
description of the partitions on the left hand side of the 
Rogers-Ramanujan identities.  This other generalization of Dyson's
rank was kindly brought to our attention by George Andrews. 

\subsection{Bressoud and Zeilberger}

A list of work connected to this proof is not complete without mentioning the bijective
Rogers-Ramanujan proof of Bressoud and Zeilberger.   In~\cite{BZ:short, BZ:long}, they give a bijection
proving Andrews' generalization of the \RR identities~(\ref{eqn:genRR}) based on the involution
principle and Bressoud's short \RR proof~\cite{Bres:easy}.  One of their maps,~$\it{\Phi}$
in~\cite{BZ:long}, acts similarly to our maps~$\mathfrak{D}^{k,m}_r$ for certain $k$, $m$, and $r$.  
Unfortunately, due to the complexity of their proofs we do not give a formal connection.   The fact
that these maps are somewhat similar does however have consequences for the question in the last
section.

\subsection{Garvan and Berkovich}

Garvan has also defined a generalized notion of rank for partitions with multiple Durfee 
squares~\cite{Gar:rank}.
Though different from our definition, his rank leads to the same generating function for partitions with 
rank at most $-r$ as we derived the previous section.  
Based on this generating function, in~\cite{BG:obs}, Berkovich and Garvan 
ask for a symmetry similar to Dyson's ``new symmetry'' for Garvan's generalized rank and for a Dyson-style proof of their 
generating function.  They note that it ``turned out to be very difficult to prove in a combinatorial fashion.''

We will explain the relationship between our generalization of rank and Garvan's definition, 
and the two symmetries associated with both definitions.  We will also be able to show why the 
Dyson-style proof sought by Berkovich and Garvan turns out to be difficult to find.

\subsubsection{Garvan's rank and conjugation}
\begin{figure}[hbt]
\begin{center}
\psfrag{2}{$\la$}
\epsfig{file=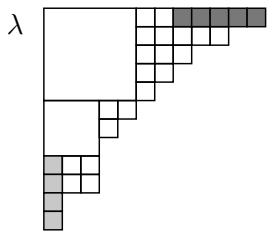,height=3cm}
\end{center}
\caption{Partition $\la = (12,10,8,7,6,5,4,3,3,3,1,1)$ has $ga_2(\la) = 5$, $gb_2(\la) = 4$, and $gr_2(\la) = 5-4=1$.}
\label{figGarRank}
\end{figure}

Recall that for a partition $\la$ with $k$ successive Durfee squares, we denote the partitions to the right of these Durfee squares by $\la^1, \la^2, ..., \la^k$ and the partition below the $k$th Durfee square by $\al$.  
\begin{defn}[Garvan,~\cite{Gar:rank}]
Let $\la$ be a partition with at least $k$ successive Durfee squares, where the $k$th Durfee square has size $N_k$.  Define
$$
\aligned
ga_k(\la) = 	& \text{ the number of columns of ${\la^1}$ whose length $\leq N_k$} \,,  \text{ and} \\
gb_k(\la) = 	& \ \ell(\al) \, .
\endaligned
$$
Also, define 
$$
gr_k(\la) = ga_k(\la) - gb_k(\la) \, .
$$
We will call $gr_k(\la)$ Garvan's rank.  
\end{defn}
Garvan called $gr_k(\la)$ the $(k+1)$-rank of $\la$.  See Figure~\ref{figGarRank} for an example.  

Garvan also described a very natural conjugation for partitions with $k$ successive Durfee squares.  
For any partition $\la$ with $k$ successive Durfee squares of size $N_1, N_2, ..., N_k$, let
\begin{itemize}
\item
$\al$ be the partition below the $k$th Durfee square, and
\item
$\be^\pr$ be the partition consisting of columns sitting to the right of the first Durfee square of $\la$ whose length is $\leq N_k$.
\end{itemize}

\begin{figure}[hbt]
\begin{center}
\psfrag{1}{$\la$}
\psfrag{2}{$\mu$}
\psfrag{a}{$\be^\pr$}
\psfrag{b}{$\al$}
\psfrag{c}{$\al^\pr$}
\psfrag{d}{$\be$}
\epsfig{file=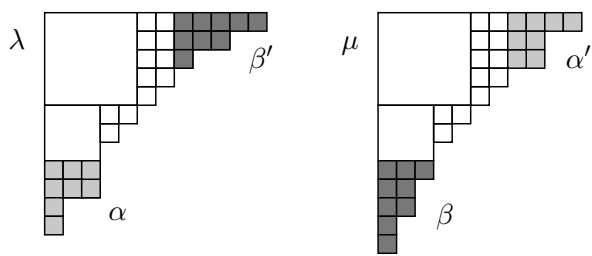,height=3.25cm}
\end{center}
\caption{Partition $\la = (12,10,8,7,6,5,4,3,3,3,1,1)$ and its Garvan conjugate $\mu = (11,9,9,7,6,5,4,3,3,2,2,1,1)$.}
\label{figGarCong}
\end{figure}

The conjugate is obtained by replacing $\al$ and $\be^\pr$ by $\be$ and $\al^\pr$, respectively.  
Note that conjugation is clearly an involution that	 sends Garvan's rank of a partition to its negative.  See Figure~\ref{figGarCong}.

\subsubsection{Relationship between $(k,m)$-rank and Garvan's ranks}

One first theorem tells us that $(k,0)$-rank and Garvan's ranks have the same 
distribution on partitions of $n$.

\begin{theorem}
\label{thm:garvrank}
For any  $k, n, s, t \in \zz$, the number of partitions $\la$ of $n$ with~$k$ successive Durfee squares, of size $N_1, N_2, ...,
N_k$, with $a_{k,0}(\la) = s$ and $b_{k,0}(\la) = t$ is equal to the number of partitions $\mu$ of $n$ with $k$ successive Durfee
squares, of size $N_1, N_2, ..., N_k$, with $ga_k(\mu) = s$ and $gb_k(\mu) = t$.  
In particular, for any  $k,n,r \in \zz$, the number of partitions $\la$ of $n$ with $r_{k,0}(\la) = r$ is equal to the number
of partitions $\mu$ of $n$ with $gr_k(\mu) = r$.  
\end{theorem}

We will give an outline of the proof of this theorem to show the relationship between the two ranks;
the full proof of this theorem is found in~\cite{Bou:Thesis}.  

\begin{figure}[hbt]
\begin{center}
\psfrag{1}{$\la$}
\psfrag{2}{$\mu$}
\psfrag{3}{$\al$}
\psfrag{4}{$\be^\pr$}
\epsfig{file=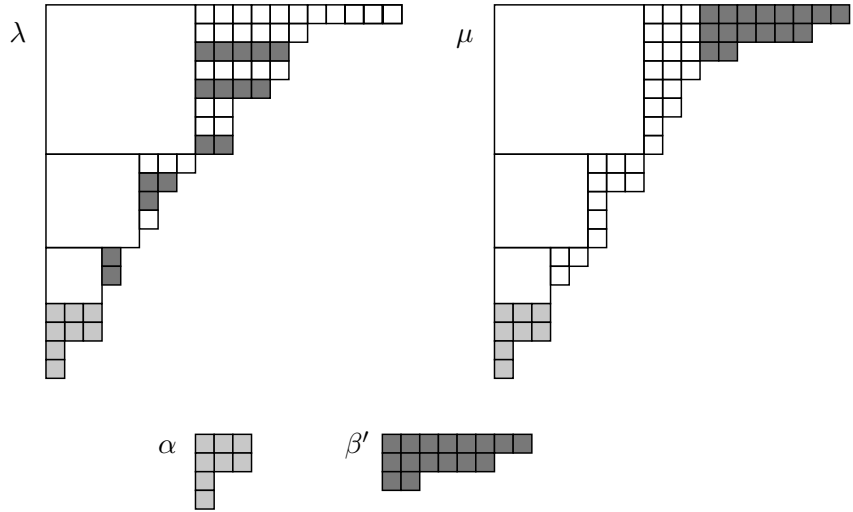,height=6.75cm}
\end{center}
\caption{The shaded squares are those exchanged by conjugation.}
\label{figJKab}
\end{figure}

For both definitions, consider the squares of the Young diagram that are exchanged by the corresponding conjugation and
the squares that are {\em not} exchanged by the corresponding conjugation.  
See Figure~\ref{figJKab}.  

The squares of $\mu$ that are exchanged by Garvan's conjugation
make up $\al$, the partition below $k$th Durfee square, and $\be^\pr$, 
the partition made up of columns having height less than or equal
to $N_k$ in~$\la^1$, the partition to the right of the first Durfee square.  
In our conguation, $\mathfrak{C}^k$, the squares of $\la$ 
that are exchanged are (in the notation from section~\ref{sec:Sym1}) $\al$, 
the partition below $k$th Durfee square, and
the partition $\be^\pr$ obtained by iteratively removing selected 
parts from $\la^1, \la^2, ..., \la^k$.  Recall that this $\be^\pr$ has
at most $N_k$ parts just as the $\be^\pr$ from Garvan's conjugation.  

\begin{figure}[hbt]
\begin{center}
\psfrag{1}{$\widetilde{\la}$}
\psfrag{2}{$\widetilde{\mu}$}
\epsfig{file=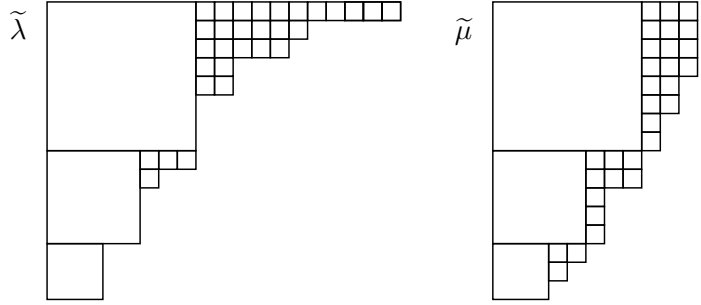,height=4cm}
\end{center}
\caption{Shaded parts have been removed to leave only the parts untouched by conjugation.}
\label{figJKEg}
\end{figure}

Now consider the squares that are {\em not} exchanged by the corresponding conjugation.  
We are left with two different sets of
partitions.  In the case of Garvan's rank, the partitions have no part below the $k$th 
Durfee square and no column whose 
length is less than or equal to $N_k$ to the right of the first Durfee square.  
In the case of our rank, the partitions have no part below the $k$th Durfee square and no 
part to the right of the
bottom $N_k$ rows of each Durfee square.  See Figure~\ref{figJKEg} for an example of these 
types of partitions.  

From this point, the proof may proceed in two ways: bijectively or by considering the generating function keeping track of $ga_k$, $gb_k$,
$gr_k$ and size of the partitions and $a_{k,0}$, $b_{k,0}$, $r_{k,0}$ and size of the partitions depending on the case.

In both cases, the squares that are exchanged by conjugation match up, $\al$ with $\al$ and $\be^\pr$ with~$\be^\pr$.   The squares
that are not exchanged by conjugation also correspond.  Notice that these latter squares do not contribute to any of the
statistics other than size of the partitions.  Therefore we only need to show that there is the same number of either type of a
given size $n$ with Durfee squares of size $N_1, ..., N_k$.  By generating functions this is a simple exercise involving
$q$-binomials.  Alternatively, a bijection can be found based a map due to Franklin, from section 20 of~\cite{Syl:cons}.

This setup also shows us that Garvan's conjugation and our conjugation are related in a natural way since we
match up exactly the squares that are exchanged by either conjugation.  

\subsubsection{A version of Dyson's map for Garvan's rank}  
Since Garvan's rank and $(k,0)$-rank are
equidistributed on partitions of a given size and 
since we can relate the conjugations associated with each rank,
it is natural look for a Dyson-like map
similar to ask our second symmetry, $\mathfrak{D}^{k,m}_r$, for Garvan's rank.  

Of course, it is possible to map $\mathfrak{D}^{k,m}_r$ through the bijection refered to in the previous section and get some
(akwardly described) map that deals with Garvan's rank.  However a few subtleties arise.   Our bijection $\mathfrak{D}^{k,m}_r$
involves both the $(k,m)$-rank of a partition and the $(k,m+2)$-rank of a partition.  Moreover, it changes the Durfee rectangle
structure of the partition.  To give a version of Dyson's map for Garvan's notion of rank based on our work neccessarily involves a
more general notion of Garvan's rank for partitions with  successive Durfee rectangles (as opposed to only squares) and such a
Dyson-like bijection neccessarily changes the Durfee rectangle structure of the partition.  

In fact, this last point is the reason that, as Berkovich and Garvan stated, it is difficult to find a Dyson's map for
Garvan's rank.  It turns out that even on fairly small examples numerical evidence shows that no Dyson-like bijection for Garvan's rank can maintain the
Durfee square structure of the partitions.

\section{Further question} 
\label{secQues}

Andrews generalized the Rogers-Ramanujan identities further than the identities~(\ref{eqn:genRR}) we have proved here.  
For $1 \leq a \leq k$, he proved that: 
\smallskip
\begin{equation}
\label{eqn:RRfurther}
\sum_{n_1 = 0}^\infty \cdots \sum_{n_{k-1} = 0}^\infty 
\frac{q^{N_1^2+N_2^2+\cdots+N_{k-1}^2 + N_a + \cdots + N_{k-1}}}{(q)_{n_1}(q)_{n_2}\cdots(q)_{n_{k-1}}} \ =
\hspace{-3mm}
\prod_{
\scriptsize{\begin{array}{c}
n=1 \\
n \not\equiv 0, \pm a~(\mathrm{mod~} 2k+1)
\end{array}}}^{\infty} \hspace{-4mm} \frac{1}{1-q^n}\end{equation}
and provided a combinatorial interpretation of the left hand side as a Durfee dissection using both Durfee squares and $1$-Durfee
rectangles~\cite{Andr:odd, Andr:Durfee}.  Further generalizations that lend themselves to similar interpretations have been given
by Bressoud as well~\cite{Bres:gen, Bres:int} and by Garrett, Ismail, and Stanton~\cite{GIS:var}.
Can our proof be extended to prove these identities?

It is fairly simple to extend our definition of $(k,m)$-rank and obtain bijections proving a first and second symmetry for these
partitions in Andrews' identity.  However, in this case, the second symmetry is not enough to determine the generating function
for partitions with rank at most $-r$.  In order to complete the proof, a new idea is required.  

On the other hand, there is evidence that our proof will not extend.  The \RR bijection given by Bressoud and
Zeilberger~\cite{BZ:short, BZ:long} is a combinatorialization of a short proof of Bressoud~\cite{Bres:easy} in which he proves
the following generalization of Schur's identity: 
\begin{equation}
\label{eqn:RRBres}
\aligned
&\sum_{s_1 = 0}^\infty \cdots \sum_{s_k = 0}^\infty 
\frac{q^{s_1^2+s_2^2+\dots+s_k^2}}{(q)_{N-s_1}(q)_{s_1-s_2}\dots(q)_{s_{k-1}-s_k}(q)_{2s_k}} (-xq;q)_{s_k} (-x^{-1};q)_{s_k} \\
&\qquad = \frac{1}{(q)_{2N}} \sum_{j = -\infty}^{\infty} x^jq^{\frac{(2k+1)j^2+j}{2}} \qbinom{2N}{N-m} \, . 
\endaligned
\end{equation}
However, this generalization is quite different from Andrews' generalization given above.  Since our map $\mathfrak{D}^{k,m}_r$
acts similarly to one of Bressoud and Zeilberger's maps, it may be that our proof is more likely to extend to this generalization
rather than equation~(\ref{eqn:RRfurther}).  

If our proof were extended to either case, this would also give a proof of the second \RR identity.  

\section{Acknowledgments} 
The author is grateful to
George Andrews, Igor Pak, and Richard Stanley for their support, encouragement, and helpful comments.

\bibliographystyle{plain}
% amsplain is also an okay style
\bibliography{DRR}

\end{document}